# EXPONENTIAL ASYMPTOTICS AND LAW OF THE ITERATED LOGARITHM FOR INTERSECTION LOCAL TIMES OF RANDOM WALKS[1]

By Xia Chen

*University of Tennessee*

Let $\alpha([0,1]^p)$ denote the intersection local time of $p$ independent $d$-dimensional Brownian motions running up to the time 1. Under the conditions $p(d-2)<d$ and $d\ge 2$, we prove

$$\lim_{t\to\infty} t^{-1}\log \mathbb{P}\{\alpha([0,1]^p)\ge t^{(d(p-1))/2}\} = -\gamma_\alpha(d,p)$$

with the right-hand side being identified in terms of the the best constant of the Gagliardo–Nirenberg inequality. Within the scale of moderate deviations, we also establish the precise tail asymptotics for the intersection local time

$$I_n = \#\{(k_1,\dots,k_p)\in[1,n]^p; S_1(k_1)=\cdots=S_p(k_p)\}$$

run by the independent, symmetric, $\mathbb{Z}^d$-valued random walks $S_1(n),\dots,S_p(n)$. Our results apply to the law of the iterated logarithm. Our approach is based on Feynman–Kac type large deviation, time exponentiation, moment computation and some technologies along the lines of probability in Banach space. As an interesting coproduct, we obtain the inequality

$$(\mathbb{E}I^m_{n_1+\cdots+n_a})^{1/p} \le \sum_{\substack{k_1+\cdots+k_a=m\\ k_1,\dots,k_a\ge 0}} \frac{m!}{k_1!\cdots k_a!}(\mathbb{E}I^{k_1}_{n_1})^{1/p}\cdots(\mathbb{E}I^{k_a}_{n_a})^{1/p}$$

in the case of random walks.

**1. Introduction.** Let $\{S_1(n)\},\dots,\{S_p(n)\}$ be symmetric independent $d$-dimensional lattice-valued random walks with the same distribution. Throughout we assume that $\{S_1(n)\},\dots,\{S_p(n)\}$ have finite second moment and that the smallest group that supports these random walks is $\mathbb{Z}^d$. Write $\Gamma$ for their covariance matrix. It is known that the trajectories of $\{S_1(n)\},\dots,\{S_p(n)\}$ intersect infinitely often if and only if $p(d-2)\le d$. In this case the intersec-

Received March 2003; revised October 2003.
[1]Supported in part by NSF Grant DMS-01-02238.
*AMS 2000 subject classifications.* 60B12, 60F10, 60F15, 60G50, 60J55, 60J65.
*Key words and phrases.* Intersection local time, large (moderate) deviations, law of the iterated logarithm, Gagliardo–Nirenberg inequality.

tion local time defined by

$$I_n = \sum_{k_1,\dots,k_p=1}^{n} \mathbb{1}_{\{S_1(k_1)=\cdots=S_p(k_p)\}}$$
$$= \#\{(k_1,\dots,k_p)\in[1,n]^p; S_1(k_1)=\cdots=S_p(k_p)\}$$

tends to $\infty$ almost surely as $n\to\infty$. An important problem is to understand the long term behavior of $\{I_n\}$ as $p(d-2)\le d$. The weak laws for $\{I_n\}$ were obtained by Le Gall (1986a, b) and Rosen (1990) [see (1.4) below for the noncritical case and Le Gall (1986b) for the critical case]. Our concern in this work is large (moderate) deviations for intersection local times and its application to the law of the iterated logarithm. The critical case defined by $p(d-2)=d$ was studied by Marcus and Rosen (1997) and Rosen (1997), and the case $d=1$ was studied by Chen and Li (2004). In this work we consider the cases defined by

$$p(d-2)<d \quad \text{and} \quad d\ge 2. \tag{1.1}$$

In other words, we consider the case $d=2$, $p\ge 2$ and the case $d=3$, $p=2$.

Another closely related object is the intersection local time generated by Brownian motions. As $d=1$, the intersection local times can be represented in terms of Brownian local times. Given a Brownian local time $L(t,x)$ and its independent copies $L_1(t,x), L_2(t,x),\dots$, the intersection local times

$$\int_{-\infty}^{\infty} L^p(t,x)\,dx \quad \text{and} \quad \int_{-\infty}^{\infty}\prod_{j=1}^{p} L_j(t,x)\,dx$$

measure the time (up to $t$) spent for self-intersection and inter-path intersection, respectively. Chen and Li (2004) observed that these two types of intersection local times have similar tail behaviors as $d=1$. This phenomenon allows us to study the mixed type of intersection local time:

$$\int_{-\infty}^{\infty}\prod_{j=1}^{m} L_j^p(t,x)\,dx.$$

Indeed, it has been established by Chen and Li (2004) that

$$\lim_{\lambda\to\infty} t^{-1}\log\mathbb{P}\left\{\int_{-\infty}^{\infty}\prod_{j=1}^{m} L_j^p(1,x)\,dx \ge t^{(mp-1)/2}\right\}$$
$$= -\frac{m}{4(mp-1)}\left(\frac{mp+1}{2}\right)^{(3-mp)/(mp-1)} B\left(\frac{1}{mp-1},\frac{1}{2}\right)^2,$$

where $B(\cdot,\cdot)$ is the beta function.

In the multidimensional case, the self-intersection local time and the inter-path intersection local time have different asymptotic behaviors and therefore are treated separately. In this paper, we deal only with the intersection local times run by $p$ independent Brownian motions. The interested reader is referred to Bass and Chen (2004) for a recent development on large deviations and the laws of the iterated logarithm for self-intersection local times of two-dimensional Brownian motions.

Let $W_1(t), \dots, W_p(t)$ be independent $d$-dimensional Brownian motions. According to Dvoretzky, Erdös and Kakutani (1950, 1954), the set of intersections

$$\bigcap_{j=1}^{p} \{x \in \mathbb{R}^d; x = W_j(t) \text{ for some } t \geq 0\}$$

contains points different from 0 if and only if $p(d-2) < d$. The interested reader is also referred to a survey paper by Khoshnevisan (2003) for an elementary proof of this result and for an overview of the whole area of path intersection.

Under condition (1.1), the intersection local time $\alpha(ds_1, \dots, ds_p)$ of $W_1(t), \dots, W_p(t)$ is defined as a random measure on $(\mathbb{R}^+)^p$ supported on

$$\{(t_1, \dots, t_p) \in (\mathbb{R}^+)^p; W_1(t_1) = \cdots = W_p(t_p)\}.$$

It is formally written as

$$\alpha(ds_1, \dots, ds_p) = \delta_0(W_1(s_1) - W_2(s_2)) \cdots \delta_0(W_{p-1}(s_{p-1}) - W_p(s_p))\, ds_1 \cdots ds_p$$

or

$$\alpha(ds_1, \dots, ds_p) = \left[\int_{\mathbb{R}^d} \prod_{j=1}^{p} \delta_0(W_j(s_j) - x)\, dx\right] ds_1 \cdots ds_p.$$

There are two equivalent ways to construct Brownian intersection local time in the multidimensional case. Geman, Horowitz and Rosen (1984) proved that under (1.1), the occupation measure on $\mathbb{R}^{d(p-1)}$ given by

$$f \mapsto \int_A f(W_1(t_1) - W_2(t_2), \dots, W_{p-1}(t_{p-1}) - W_p(t_p))\, dt_1 \cdots dt_p$$

is absolutely continuous with respect to Lebesgue measure on $\mathbb{R}^{d(p-1)}$ for any Borel set $A \subset (\mathbb{R}^p)^+$ and the density $\alpha(x, A)$ of the occupation measure can be chosen in such a way that the function

$$(x, t_1, \dots, t_p) \mapsto \alpha(x, [0, t_1] \times \cdots \times [0, t_p]),$$
$$x \in \mathbb{R}^{d(p-1)}, (t_1, \dots, t_p) \in (\mathbb{R}^p)^+,$$

is jointly continuous [see Bass and Khoshnevisan (1993) for a further discussion on Hölder continuity of the Brownian intersection local times]. The random measure $\alpha(\cdot)$ on $(\mathbb{R}^p)^+$ is defined as

$$\alpha(A)=\alpha(0,A) \qquad \forall \text{Borel set } A\subset(\mathbb{R}^p)^+. \tag{1.2}$$

Another approach involves the approximation of the Dirac function. Given $\varepsilon>0$, define the random measure $\alpha_\varepsilon(\cdot)$ on $(\mathbb{R}^p)^+$ by

$$\alpha_\varepsilon(ds_1,\dots,ds_p)=\left[\int_{\mathbf{R}^d}\prod_{j=1}^p\frac{1}{C_d\varepsilon^d}\mathbb{1}_{\{|W_j(s_j)-x|\le\varepsilon\}}\,dx\right]ds_1\cdots ds_p,$$

where $C_d$ is the volume of the $d$-dimensional unit ball. When $d=2$ and $p\ge 2$, Le Gall [(1990), Theorem 1, page 183] showed that there is a unique random measure $\alpha(\cdot)$ on $(\mathbb{R}^p)^+$ such that

$$\lim_{\varepsilon\to 0^+}\alpha_\varepsilon(A_1\times\cdots\times A_p)=\alpha(A_1\times\cdots\times A_p)$$

holds in $L^m$ norm for any $m\ge 1$ and for any Borel box $A_1\times\cdots\times A_p\subset(\mathbb{R}^p)^+$. By closely examining his argument, we can see that the same conclusion applies to all cases defined by (1.1) without extra work.

By the scaling property of Brownian motions,

$$\alpha([0,t]^p)\overset{d}{=}t^{(2p-d(p-1))/2}\alpha([0,1]^p). \tag{1.3}$$

Very similarly,

$$n^{-(2p-d(p-1))/2}I_n\xrightarrow{d}\det(\Gamma)^{-(p-1)/2}\alpha([0,1]^p). \tag{1.4}$$

Indeed, it is known [see, e.g., Le Gall (1986a) and Rosen (1990)] that

$$n^{-(2p-d(p-1))/2}I_n\xrightarrow{d}\alpha_U([0,1]^p),$$

where $\alpha_U([0,t]^p)$ denotes the analogue of $\alpha([0,t]^p)$ with $W_1(t),\dots,W_p(t)$ replaced by the independent $d$-dimensional Lévy Gaussian processes $U_1(t),\dots,U_p(t)$ whose covariance matrix is $\Gamma$. By the fact that $W\overset{d}{=}\Gamma^{-1/2}U$, we have

$$\{\alpha_U([0,t]^p);t\ge 0\}\overset{d}{=}\{\det(\Gamma)^{-(p-1)/2}\alpha([0,t]^p);t\ge 0\}, \tag{1.5}$$

from which (1.4) follows. Here we recommend the interested reader to Bass and Khoshnevisan (1992) for further discussion on the uniformity (over time and initial points) of the convergence given in (1.4). So it is expected that $I_n$ and $\alpha([0,t]^p)$ have similar long term behavior.

Under the condition $p(d-2)<d$, König and Mörters (2002) established the large deviation principle for projected intersection local times,

$$\lim_{t\to\infty}t^{-1/p}\log\mathbb{P}\{l(U)\ge t\}=-C(p,d,U), \tag{1.6}$$

where

$$l(U) = \int_U dy \prod_{j=1}^{p} \int_0^{T_j} ds\, \delta_y(W_j(s)),$$

where $T_1, \ldots, T_p$ are the exit times of, respectively, the Brownian motions $W_1, \ldots, W_p$ from a ball $B(0, R)$ of center 0 and radius $R > 0$, where $U \subset B(0, R)$ is a bounded open subset that contains starting points of $W_1, \ldots, W_p$ and where the constant $C(p, d, U) > 0$ is given in terms of certain variational problems. Our situation is different from that studied by König and Mörters (2002): We let the Brownian motion run up to a deterministic time rather than random time, and our intersection local time is defined on the whole space $\mathbb{R}^d$, while theirs is limited to a bounded domain. On the other hand, we point to (Remark 4.1) a possible connection between our work and that by König and Mörters.

To see what we expect, we recall some work done by Le Gall (1986a). In Lemma 2.2 of Le Gall (1986a) it is proved that as $d = 2$ and $p \geq 2$,

$$C_1^m (m!)^{p-1} \leq \mathbb{E}\alpha([0,1]^p)^m \leq C_2^m (m!)^{p-1} (\log m)^m, \tag{1.7}$$

and as $d = 3$ and $p = 2$,

$$C_1^m (m!)^{3/2} \leq \mathbb{E}\alpha([0,1]^2)^m \leq C_2^m (m!)^{3/2}, \tag{1.8}$$

where the constants $C_1 > 0$ and $C_2 > 0$ depend only on $(d, p)$. To improve the upper bound in (1.7), we propose the following simple treatment, which works for all cases defined by (1.1).

Let $\tau_1, \ldots, \tau_p$ be i.i.d. exponential times with parameter 1, and assume the independence between $\{\tau_1, \ldots, \tau_p\}$ and $\{X_1(t), \ldots, X_p(t)\}$. Let $\Sigma_m$ be the set of the permutations on $\{1, \ldots, m\}$ and let $p_t(x)$ $(t \geq 0)$ be the $d$-dimensional Brownian densities. By Le Gall's moment formula [(1), page 182 in Le Gall (1990)],

$$\begin{aligned}
&\mathbb{E}[\alpha([0,\tau_1] \times \cdots \times [0,\tau_p])^m] \\
&\quad = \int_{(\mathbb{R}^d)^m} dx_1 \cdots dx_m \Bigg[\int_0^\infty e^{-t}\, dt \\
&\qquad\qquad \times \int_{0 \leq s_1 < \cdots < s_m \leq t} \sum_{\sigma \in \Sigma_m} \prod_{k=1}^m p_{s_k - s_{k-1}}(x_{\sigma(k)} - x_{\sigma(k-1)})\Bigg]^p \\
&\quad = \int_{(\mathbb{R}^d)^m} dx_1 \cdots dx_m \Bigg[\sum_{\sigma \in \Sigma_m} \prod_{k=1}^m \int_0^\infty e^{-t} p_t(x_{\sigma(k)} - x_{\sigma(k-1)})\, dt\Bigg]^p \\
&\quad \leq (m!)^p \int_{(\mathbb{R}^d)^m} dx_1 \cdots dx_m \Bigg[\prod_{k=1}^m \int_0^\infty e^{-t} p_t(x_k)\, dt\Bigg]^p
\end{aligned}$$

$$= (m!)^p \left[\int_{\mathbb{R}^d} \left(\int_0^\infty e^{-t} p_t(x)\, dt\right)^p dx\right]^m,$$

where the convention that $x_{\sigma(0)} = 0$ and $s_0 = 0$ is adapted, where the inequality follows from Hölder's inequality together with some suitable variable substitutions and where the second equality follows from the fact that

$$\begin{aligned} &\int_0^\infty e^{-t}\, dt \int_{0\le s_1<\cdots<s_n\le t} ds_1\cdots ds_p \varphi_1(s_1)\prod_{k=2}^m \varphi_k(s_k - s_{k-1}) \\ &\quad = \prod_{k=1}^m \int_0^\infty e^{-t}\varphi_k(t)\, dt. \end{aligned} \tag{1.9}$$

We give a short proof here. By the substitution

$$t_1 = s_1, \qquad t_2 = s_2 - s_1, \ldots, t_m = s_m - s_{m-1} \quad \text{and} \quad t_{m+1} = t - s_m,$$

the integral on the left-hand side is equal to

$$\begin{aligned} &\int_0^\infty \cdots \int_0^\infty dt_1\cdots dt_{m+1} \exp(-t_{m+1}) \prod_{k=1}^m \exp(-t_k)\varphi_k(t_k) \\ &\quad = \prod_{k=1}^m \int_0^\infty e^{-t}\varphi_k(t)\, dt. \end{aligned}$$

Notice that

$$\begin{aligned} &\int_{\mathbb{R}^d}\left(\int_0^\infty e^{-t}p_t(x)\, dt\right)^p dx \\ &\quad = \int_0^\infty\cdots\int_0^\infty dt_1\cdots dt_p \exp(-(t_1+\cdots+t_p)) \int_{\mathbb{R}^d}\prod_{j=1}^p p_{t_j}(x)\, dx \\ &\quad = (2\pi)^{-(d(p-1))/2} \\ &\qquad \times \int_0^\infty\cdots\int_0^\infty \exp(-(t_1+\cdots+t_p)) \left(\sum_{j=1}^p \prod_{\substack{1\le k\le p \\ k\ne j}} t_k\right)^{-d/2} dt_1\cdots dt_p. \end{aligned}$$

By arithmetic–geometric mean inequality,

$$\frac{1}{p}\sum_{j=1}^p \prod_{\substack{1\le k\le p \\ k\ne j}} t_k \ge \sqrt[p]{\prod_{j=1}^p \prod_{\substack{1\le k\le p \\ k\ne j}} t_k} = \prod_{j=1}^p t_j^{(p-1)/p}.$$

Notice that $d(p-1)/(2p) < 1$ under (1.1). So we have

$$\int_{\mathbb{R}^d}\left(\int_0^\infty e^{-t}p_t(x)\, dt\right)^p dx$$

$$\le (2\pi)^{-(d(p-1))/2} p^{-d/2}\left(\int_0^\infty t^{-(d(p-1))/(2p)} e^{-t}\,dt\right)^p < \infty.$$

On the other hand, $\tau_{\min} \equiv \min\{\tau_1, \dots, \tau_p\}$ is an exponential time with parameter $p$. By the scaling property given in (1.3),

$$\begin{aligned}
&\mathbb{E}[\alpha([0,\tau_1]\times\cdots\times[0,\tau_p])^m] \\
&\quad\ge \mathbb{E}[\alpha([0,\tau_{\min}]^p)^m] = \mathbb{E}\tau_{\min}^{((2p-d(p-1))/2)m}\mathbb{E}[\alpha([0,1]^p)^m] \qquad (1.10)\\
&\quad= p^{-(((2p-d(p-1))/2)m-1)}\Gamma\left(1+\frac{2p-d(p-1)}{2}m\right)\mathbb{E}\alpha([0,1]^p)^m.
\end{aligned}$$

Summarizing what we have,

$$\begin{aligned}
\mathbb{E}\alpha([0,1]^p)^m &\le p^{(((2p-d(p-1))/2)m+1)}\Gamma\left(1+\frac{2p-d(p-1)}{2}m\right)^{-1}\\
&\qquad\times (m!)^p\left[\int_{\mathbb{R}^d}\left(\int_0^\infty e^{-t}p_t(x)\,dt\right)^p dx\right]^m.
\end{aligned}$$

By Stirling's formula,

$$\begin{aligned}
&\limsup_{m\to\infty}\sqrt[m]{\frac{\mathbb{E}\alpha([0,1]^p)^m}{(m!)^{(d(p-1))/2}}}\\
&\quad\le\left(\frac{2p-d(p-1)}{2p}\right)^{-(2p-d(p-1))/2}\int_{\mathbb{R}^d}\left(\int_0^\infty e^{-t}p_t(x)\,dt\right)^p dx,
\end{aligned}$$

$$\begin{aligned}
&\limsup_{k\to\infty}\sqrt[k]{\frac{\mathbb{E}\alpha([0,1]^p)^{2k/(d(p-1))}}{k!}} \qquad (1.11)\\
&\quad\le\left(\frac{2p-d(p-1)}{2p}\right)^{-(2p-d(p-1))/(d(p-1))}\\
&\qquad\times\left(\frac{d(p-1)}{2}\right)^{-1}\left(\int_{\mathbb{R}^d}\left(\int_0^\infty e^{-t}p_t(x)\,dt\right)^p dx\right)^{2/(d(p-1))}.
\end{aligned}$$

By the first estimate, we obtain an upper bound for $\mathbb{E}\alpha([0,1]^p)^m$, which claims, together with the lower bounds given in (1.7) and (1.8), that

$$C_1^m(m!)^{(d(p-1))/2} \le \mathbb{E}\alpha([0,1]^p)^m \le C_2^m(m!)^{(d(p-1))/2}. \qquad (1.12)$$

Consequently, there is a constant $\gamma_\alpha(d,p) > 0$, such that

$$\mathbb{E}\exp\{\gamma\alpha([0,1]^p)^{2/(d(p-1))}\}\begin{cases} <\infty, & \gamma<\gamma_\alpha(d,p),\\ =\infty, & \gamma>\gamma_\alpha(d,p).\end{cases} \qquad (1.13)$$

In the special case $d = p = 2$, (1.13) was obtained by Le Gall (1994).

QUESTION. What can we say about the critical exponent $\gamma_\alpha(d,p)$? Can we strengthen (1.13) into the large deviation

$$\lim_{t\to\infty} t^{-1}\log\mathbb{P}\{\alpha([0,1]^p)\ge t^{(d(p-1))/2}\} = -\gamma_\alpha(d,p)?$$

By (1.11) and the Taylor expansion, a partial answer to the question is

$$\gamma_\alpha(d,p) \ge \frac{d(p-1)}{2}\Big(\frac{2p-d(p-1)}{2p}\Big)^{(2p-d(p-1))/(d(p-1))} \times \Big(\int_{\mathbb{R}^d}\Big(\int_0^\infty e^{-t}p_t(x)\,dt\Big)^p dx\Big)^{-2/(d(p-1))}. \tag{1.14}$$

Unfortunately, (1.14) cannot be developed into an equality. The complete answer is given in Theorem 2.1.

**2. Main results.** Recall our assumption (1.1). To identify the constants given in our main results, we consider the Gagliardo–Nirenberg type inequalities

$$\|f\|_{2p} \le C\|\nabla f\|_2^{(d(p-1))/(2p)}\cdot\|f\|_2^{1-(d(p-1))/(2p)}, \qquad f\in W^{1,2}(\mathbb{R}^d), \tag{2.1}$$

where $C>0$ is a constant independent of $f$ and for each $r\ge 1$, $W^{1,r}(\mathbb{R}^d)$ denotes the Sobolev space

$$W^{1,r}(\mathbb{R}^d) = \{f\in\mathcal{L}^r(\mathbb{R}^d); \nabla f\in\mathcal{L}^r(\mathbb{R}^d)\}.$$

The validity of (2.1) can be derived from the Sobolev inequality [see, e.g., Ziemer (1989), Theorem 2.4.1, page 56]

$$\|f\|_{r^*} \le K\|\nabla f\|_r, \qquad f\in W^{1,r}(\mathbb{R}^d),$$

where $1\le r<d$, $r^* = dr(d-r)^{-1}$ and $K>0$ is a constant that depends only on $r$ and $d$. Indeed, taking $r = d(p-1)p^{-1}$, we have

$$\|f\|_{d(p-1)} \le K\|\nabla f\|_{d(p-1)p^{-1}}. \tag{2.2}$$

Replacing $f$ with $|f|^{2p/(d(p-1))}$, we have

$$\begin{aligned}\|f\|_{2p} &\le \Big(\frac{2p}{d(p-1)}K\Big)^{(d(p-1))/(2p)} \\ &\quad\times\Big(\int_{\mathbb{R}^d}|\nabla f(x)|^{(d(p-1))/p}|f(x)|^{(2p-d(p-1))/p}\,dx\Big)^{1/2} \\ &\le \Big(\frac{2p}{d(p-1)}K\Big)^{(d(p-1))/(2p)}\|\nabla f\|_2^{(d(p-1))/(2p)}\cdot\|f\|_2^{1-(d(p-1))/(2p)},\end{aligned}$$

where the second step follows from Hölder inequality.

Let $\kappa(d,p)$ be the best constant of the Gagliardo–Nirenberg inequality given in (2.1):

$$\kappa(d,p) = \inf\{C>0; \|f\|_{2p} \le C\|\nabla f\|_2^{(d(p-1))/(2p)} \cdot \|f\|_2^{1-(d(p-1))/(2p)} \quad \text{for } \forall f \in W^{1,2}(\mathbb{R}^d)\}.$$

Then $0 < \kappa(d,p) < \infty$.

THEOREM 2.1. *Under condition* (1.1),

$$\lim_{t\to\infty} t^{-1} \log \mathbb{P}\{\alpha([0,1]^p) \ge t^{(d(p-1))/2}\} = -\frac{p}{2}\kappa(d,p)^{-4p/(d(p-1))}. \tag{2.3}$$

Finding the best constants for Gagliardo–Nirenberg inequalities appears to be a difficult problem which remains open in general. It has been attracting considerable attention partially due to its connection to some problems in physics. The best constant for Nash's inequality, which is a special case of Gagliardo–Nirenberg inequalities, was found by Carlen and Loss (1993). See also papers by Cordero-Erausquin, Nazaret and Villani (2004) and by Del Pino and Dolbeault (2002) for recent progress on the best constants for a class of Gagliardo–Nirenberg inequalities. See the paper by Del Pino and Dolbeault (2003) for a connection between the best constants for Gagliardo–Nirenberg inequalities and logarithmic Sobolev inequalities. Two papers are directly related to $\kappa(d,p)$: In the case $d=2$ and $p=3$, Levine (1980) obtained the sharp estimate

$$\sqrt[3]{\frac{1}{4.6016}} < \kappa(2,3) < \sqrt[3]{\frac{1}{4.5981}}.$$

He conjectured that $\kappa(2,3) = \pi^{-4/9}$. Weinstein (1983) studied the problem of the best constants for the Gagliardo–Nirenberg inequalities of the type (2.1). It was shown [Weinstein (1983), Theorem B] that under (1.1), the best constant is attained at an infinitely smooth, positive and spherically symmetric function $f_0$, which solves the nonlinear equation

$$\frac{d(p-1)}{2}\triangle f - \frac{2p-d(p-1)}{2} f + f^{2p-1} = 0.$$

In addition, $f_0$ has the smallest $L^2$ norm among all solutions of the above equation (such a solution is called a ground state solution). Furthermore [Weinstein (1983), (I.3)],

$$\kappa(d,p) = (p\|f_0\|_2^{-2(p-1)})^{1/(2p)}.$$

Using this result, Weinstein (1983) obtained the following numerical approximation in the case $d=p=2$:

$$\kappa(2,2) \approx \sqrt[4]{\frac{1}{\pi \times (1.86225\dots)}}.$$

By comparing Theorem 2.1 with (1.14), the following bound of $\kappa(d,p)$ is a by-product:

$$\kappa(d,p) \le \left(\frac{p}{d(p-1)}\right)^{(d(p-1))/(4p)} \left(\frac{2p}{2p-d(p-1)}\right)^{(2p-d(p-1))/(4p)} \times \left(\int_{\mathbb{R}^d} \left(\int_0^\infty e^{-t} p_t(x)\,dt\right)^p dx\right)^{1/(2p)}. \tag{2.4}$$

In Lemma A.2, we connect $\kappa(d,p)$ to a variational problem.

We now turn to the random walks. Let $\{b_n\}$ be a positive sequence that satisfies

$$b_n \to \infty \quad \text{and} \quad b_n/n \to 0 \qquad (n \to \infty). \tag{2.5}$$

Theorem 2.2. *Under conditions* (1.1) *and* (2.5),

$$\lim_{n\to\infty} \frac{1}{b_n} \log \mathbb{P}\{I_n \ge \lambda n^{(2p-d(p-1))/2} b_n^{(d(p-1))/2}\} = -\frac{p}{2} \det(\Gamma)^{1/d} \kappa(d,p)^{-4p/(d(p-1))} \lambda^{2/(d(p-1))} \tag{2.6}$$

*for all* $\lambda > 0$.

Our large and moderate deviations apply to the following law of the iterated logarithm.

Theorem 2.3. *Under condition* (1.1),

$$\limsup_{t\to\infty} t^{-(2p-d(p-1))/2} (\log\log t)^{-(d(p-1))/2} \alpha([0,t]^p) = \left(\frac{2}{p}\right)^{(d(p-1))/2} \kappa(d,p)^{2p} \qquad \textit{a.s.}, \tag{2.7}$$

$$\limsup_{n\to\infty} n^{-(2p-d(p-1))/2} (\log\log n)^{-(d(p-1))/2} I_n = \left(\frac{2}{p}\right)^{(d(p-1))/2} \det(\Gamma)^{-(p-1)/2} \kappa(d,p)^{2p} \qquad \textit{a.s.} \tag{2.8}$$

*Specifically,*

$$\limsup_{t\to\infty} \frac{1}{t(\log\log t)^{p-1}} \alpha([0,t]^p) = \left(\frac{2}{p}\right)^{p-1} \kappa(2,p)^{2p} \qquad \textit{a.s.},$$

$$\limsup_{n\to\infty} \frac{1}{n(\log\log n)^{p-1}} I_n = \left(\frac{2}{p}\right)^{p-1} \det(\Gamma)^{-(p-1)/2} \kappa(2,p)^{2p} \qquad \textit{a.s.}$$

*as* $d=2$ *and* $p\ge 2$, *and*

$$\limsup_{t\to\infty}\frac{1}{\sqrt{t(\log\log t)^3}}\alpha([0,t]^p)=\kappa(3,2)^4\qquad a.s.,$$

$$\limsup_{n\to\infty}\frac{1}{\sqrt{n(\log\log n)^3}}I_n=\det(\Gamma)^{-1/2}\kappa(3,2)^4\qquad a.s.$$

*as* $d=3$ *and* $p=2$.

Recall that the trajectories of $\{S_1(n)\},\dots,\{S_p(n)\}$ intersect infinitely often if and only if $p(d-2)\le d$. In the critical cases defined as $p(d-2)=d$—the case $d=4$, $p=2$ and the case $d=p=3$—the law of the iterated logarithm was obtained by Marcus and Rosen (1997) and Rosen (1997), respectively. Under the assumption of finite third moment, it was proved [Marcus and Rosen (1997)] that

$$\limsup_{n\to\infty}\frac{I_n}{\log n\log\log\log n}=\frac{1}{2\pi^2\sqrt{\det(\Gamma)}}\qquad \text{a.s.}\tag{2.9}$$

as $d=4$ and $p=2$, and [Rosen (1997)] that

$$\limsup_{n\to\infty}\frac{I_n}{\log n\log\log\log n}=\frac{1}{\pi\det(\Gamma)}\qquad \text{a.s.}\tag{2.10}$$

as $d=p=3$.

The case $d=1$ was studied by Chen and Li (2004). As a special form of a general result given in Theorem 1.4 of Chen and Li (2004),

$$\begin{aligned}&\limsup_{n\to\infty}n^{-(p+1)/2}(\log\log n)^{-(p-1)/2}I_n\\&\quad=\left(\frac{4(p-1)}{p\sigma^2}\right)^{(p-1)/2}\left(\frac{p+1}{2}\right)^{(p-3)/2}B\left(\frac{1}{p-1},\frac{1}{2}\right)^{-(p-1)}\qquad \text{a.s.}\end{aligned}\tag{2.11}$$

as $d=1$ and $p\ge 2$, where $\sigma^2>0$ is the variance of the random walks.

The law of the iterated logarithm given in Theorem 2.3 solves the cases left by the previous works.

In addition to being important in the study of random paths, the notion of intersection local times of independent Brownian trajectories is connected to some other interesting problems. Bass and Chen (2004) proved that the renormalized 2-multiple self-intersection local time run by a two-dimensional Brownian motion satisfies exactly the same large deviation and the law of the iterated logarithm as does $\alpha([0,t]^2)$ in the case $d=p=2$. Chen and Li (2004) pointed out how the intersection local times are related to the local times of additive Lévy processes [see Khoshnevisan, Xiao and Zhang (2003a, b) for some later developments in this area]. In their paper, König and Mörters (2002) applied the large deviation result given in (1.6) to the problem of finding the Hausdorff dimension spectrum for thick points of the intersection of two independent Brownian paths in $\mathbb{R}^3$.

The most interesting link is the range problem. While the intersection local time $I_n$ counts the times spent in intersecting, the random variable

$$J_n = \#\{S_1[1,n] \cap \cdots \cap S_p[1,n]\} \tag{2.12}$$

counts the sites of intersection. Clearly, $J_n \le I_n$ and the difference is caused by the possibility that the trajectories intersect more than once at the same site. The weak laws for $J_n$ were studied by Le Gall (1986a, b) and Le Gall and Rosen (1991), and it has been observed that in the transient case ($d \ge 3$), $J_n$ behaves like $\gamma^p I_n$, where $\gamma = \mathbb{P}\{S_n \ne 0 \ \forall n \ge 1\}$. In light of (2.9) and (2.10), therefore, the law of the iterated logarithm was given by Marcus and Rosen (1997) and Rosen (1997) as $p(d-2) = d$. On the other hand, Le Gall (1986a) proved that as $d = p = 2$,

$$\frac{(\log n)^2}{n} J_n \xrightarrow{d} (2\pi)^2 \det(\Gamma)^{-1/2} \alpha([0,1]^2). \tag{2.13}$$

Comparing (2.13) with (1.4), we observe a sharp contrast between $I_n$ and $J_n$ as $d = 2$. Theorem 2.1 established in our work opens the possibility to understand the tail behavior of $J_n$ under the condition (1.1). We leave this to future study.

We outline some key technique points in each of the following sections. In Section 3, we establish the large (moderate) deviations for the intersection local times smoothed by convolution. Different from the one-dimensional case, multidimensional Brownian motions do not have local time. Even in the case of random walks, the absence of a modulus of continuity for local times makes it difficult to handle intersection local time directly as we did before [Chen and Li (2004)]. For this reason, the large (moderate) deviations are accomplished first for the $L_p$ norms and multilinear forms of the occupation times over small balls instead of the local times. In spite of the difference mentioned above, some of the ideas developed by Chen and Li (2004) are adapted here: A Feynman–Kac type large deviation given by Remillard (2000) for Brownian occupation times and analogous minorization established by Chen and Li (2004) in the context of random walks, an idea of localization developed by Donsker and Varadhan (1975) and Mansmann (1991), a deterministic comparison (via Hölder type inequality) between multilinear form and $L_p$ norms, a way to establish exponential tightness introduced by de Acosta (1985) and $L_p$ embedding are the main ingredients in the proof of Theorem 3.1. We point out that Donsker and Varadhan's (1974) large deviation principle for empirical processes could be used [see Mansmann (1991) for "how" in the case $d = 1$] in the proof of (3.2) and (3.3). We chose not to do so because the proof of (3.4) and (3.5) demands an approach which can be extended (at least partially) to the case of random walks.

In Section 4, we prove the upper bound for Theorem 2.1. The idea is approximation via Theorems 3.1. Le Gall's moment formula and time exponentiation (Laplacian transform) are essential tools. Laplacian transform has been developed into an important tool in the study of limit laws for occupation times since the remarkable work by Darling and Kac (1957). The interested reader is referred to the survey paper by Fitzsimmons and Pitman (1999) for an overview. It is worth mention that our situation is not quite standard: We have to deal with $p$ independent exponential times at the same time.

In Section 5, we prove the upper bound for Theorem 2.2, which is by no means a trivial consequence of Theorem 2.1 and the invariance principle, due to the discontinuity of the functionals involved. The treatment we present is completely different from that for the upper bound of Theorem 2.1. The central piece of our approach is a moment inequality given in Theorem 5.1. This inequality appears to be interesting for its own sake.

In Section 6, we prove the lower bounds for both Theorems 2.1 and 2.2. For the needs of the law of the iterated logarithm, the statement given in Theorem 6.1 is more than we need for Theorems 2.1 and 2.2: We allow the independent paths to start at different points and we establish the lower bounds uniformly over the starting points. The key step is to show that the moments of intersection local times can be minorized by those of the multilinear forms that appear in Theorem 3.1.

In Section 7, we prove Theorem 2.3. With the extensive preparation on the tail estimate, the proof is a standard practice of the Borel–Cantelli lemma. In the Appendix, we give two analytic lemmas which are used in the proof of our main theorems.

**3. Smoothing intersection local times.** To simplify the notation, we denote $W(t)$ for a $d$-dimensional Brownian motion and $\{S(n)\}$ for a $d$-dimensional random walk with the same distribution as $\{S_1(n)\},\dots,\{S_p(n)\}$, whenever only a single Brownian motion or a single random walk is involved. In this section, we let the small number $\varepsilon>0$ be fixed but arbitrary and write

$$L(t,x,\varepsilon)=\frac{1}{C_d\varepsilon^d}\int_0^t \mathbb{1}_{\{|W(s)-x|\le\varepsilon\}}\,ds, \qquad x\in\mathbb{R}^d, t\ge 0,$$

$$l(n,x,\varepsilon)=\frac{1}{\#\{B_n\}}\sum_{k=1}^n \mathbb{1}_{\{S(k)-x\in B_n\}}, \qquad x\in\mathbb{Z}^d,\ n=1,2,\dots,$$

where $C_d$ is the volume of the $d$-dimensional unit ball and

$$B_n=\{y\in\mathbb{Z}^d;\ |y|\le\varepsilon\sqrt{n/b_n}\,\}.$$

For each $1 \le j \le p$, let $L_j(t,x,\varepsilon)$ and $l_j(n,x,\varepsilon)$ be the analogues of $L(t,x,\varepsilon)$ and $l(n,x,\varepsilon)$ with $W(t)$ and $S(n)$ being replaced by $W_j(t)$ and $S_j(n)$, respectively.

For any locally integrable function $f$ on $\mathbb{R}^d$, we introduce the notation

$$f_\varepsilon(x) = \frac{1}{C_d\varepsilon^d}\int_{\{|y-x|\le\varepsilon\}} f(y)\,dy.$$

Given $\theta > 0$, write

$$M_\varepsilon(\theta) = \sup_{f\in\mathcal{F}_d}\biggl\{\theta\biggl(\int_{\mathbb{R}^d}[(f^2)_\varepsilon(x)]^p\,dx\biggr)^{1/p} - \tfrac{1}{2}\int_{\mathbb{R}^d}|\nabla f(x)|^2\,dx\biggr\},$$

$$N_\varepsilon(\theta) = \sup_{f\in\mathcal{F}_d}\biggl\{\theta\biggl(\int_{\mathbb{R}^d}[(f^2)_\varepsilon(x)]^p\,dx\biggr)^{1/p} - \frac{p}{2}\int_{\mathbb{R}^d}|\nabla f(x)|^2\,dx\biggr\},$$

$$\widetilde{M}_\varepsilon(\theta) = \sup_{f\in\mathcal{F}_d}\biggl\{\theta\biggl(\int_{\mathbb{R}^d}[(f^2)_\varepsilon(x)]^p\,dx\biggr)^{1/p} - \tfrac{1}{2}\int_{\mathbb{R}^d}\langle\nabla f,\Gamma\,\nabla f\rangle\,dx\biggr\},$$

$$\widetilde{N}_\varepsilon(\theta) = \sup_{f\in\mathcal{F}_d}\biggl\{\theta\biggl(\int_{\mathbb{R}^d}[(f^2)_\varepsilon(x)]^p\,dx\biggr)^{1/p} - \frac{p}{2}\int_{\mathbb{R}^d}\langle\nabla f,\Gamma\,\nabla f\rangle\,dx\biggr\},$$

where

$$\mathcal{F}_d = \biggl\{f\in W^{1,2}(\mathbb{R}^d);\int_{\mathbb{R}^d}|f(x)|^2\,dx = 1\biggr\}. \tag{3.1}$$

By the fact that

$$\int_{\mathbb{R}^d}[(f^2)_\varepsilon(x)]^p\,dx \le \sup_{x\in\mathbb{R}^d}[(f^2)_\varepsilon(x)]^{p-1} \le \biggl(\frac{1}{C_d\varepsilon^d}\biggr)^{p-1},$$

the functions $M_\varepsilon(\cdot)$, $M'_\varepsilon(\cdot)$, $N_\varepsilon(\cdot)$ and $N'_\varepsilon(\cdot)$ are continuous for any fixed $\varepsilon > 0$.

THEOREM 3.1. *For any $\theta > 0$ and integers $d \ge 1$, $p \ge 2$,*

$$\lim_{t\to\infty}\frac{1}{t}\log\mathbb{E}\exp\biggl\{\theta\biggl(\int_{\mathbb{R}^d}L^p(t,x,\varepsilon)\,dx\biggr)^{1/p}\biggr\} = M_\varepsilon(\theta), \tag{3.2}$$

$$\lim_{t\to\infty}\frac{1}{t}\log\mathbb{E}\exp\biggl\{\theta\biggl(\int_{\mathbb{R}^d}\prod_{j=1}^p L_j(t,x,\varepsilon)\,dx\biggr)^{1/p}\biggr\} = N_\varepsilon(\theta), \tag{3.3}$$

$$\begin{aligned}\lim_{n\to\infty}\frac{1}{b_n}\log\mathbb{E}\exp\biggl\{\theta\biggl(\frac{b_n}{n}\biggr)^{(2p-d(p-1))/(2p)}\\ \times\biggl(\sum_{x\in\mathbb{Z}^d}l^p(n,x,\varepsilon)\biggr)^{1/p}\biggr\} = \widetilde{M}_\varepsilon(\theta),\end{aligned} \tag{3.4}$$

$$\lim_{n\to\infty}\frac{1}{b_n}\log\mathbb{E}\exp\biggl\{\theta\biggl(\frac{b_n}{n}\biggr)^{(2p-d(p-1))/(2p)} \times\biggl(\sum_{x\in\mathbb{Z}^d}\prod_{j=1}^{p} l_j(n,x,\varepsilon)\biggr)^{1/p}\biggr\} = \widetilde{N}_\varepsilon(\theta). \tag{3.5}$$

Remark 3.1. It should be emphasized that Theorem 3.1 holds for all $d \ge 1$ and all integers $p \ge 2$ [in other words, condition (1.1) plays no role here], and that, smoothed by the uniform distribution over a small ball, the self-intersection and interpath intersection present almost the same behavior. This is quite contrary to the strong dimension dependence of the intersection local times and contrary to the substantial difference in asymptotic behaviors between self-intersection local times and interpath intersection local times as $d \ge 2$. Here is our explanation: First, replacing the trajectories by their "$\varepsilon$ sausage" makes intersection always possible regardless of the values of $d$ and $p$. Second, the behavior of the intersection local times is determined by the degree of their singularity, which depends on a combination of the space structure (dimension $d$) and the pattern of intersection. The smoothing procedure eliminates the singularity and therefore eliminates the difference in behavior. Furthermore, we can see from the proof below that $p$ can be any real number larger than 1 in (3.2) and (3.4).

Proof of Theorem 3.1. We first deal with Brownian case. We start with a result based on the Feynman–Kac formula [see, e.g., Remillard (2000)],

$$\lim_{t\to\infty}\frac{1}{t}\log\mathbb{E}\exp\biggl\{\int_0^t f(W(s))\,ds\biggr\} = \sup_{g\in\mathcal{F}_d}\biggl\{\int_{\mathbb{R}^d} f(x)g^2(x)\,dx - \frac{1}{2}\int_{\mathbb{R}^d}|\nabla g(x)|^2\,dx\biggr\}, \tag{3.6}$$

where $f$ can be any bounded, measurable function $f$ on $\mathbb{R}^d$.

We now prove the upper bound for (3.2). We may let $t \to \infty$ along the integer points in our argument when it is needed. The basic idea is to view $\{L(t,\cdot,\varepsilon); t \ge 0\}$ as a process taking values in the Banach space $\mathcal{L}^p(\mathbb{R}^d)$. Then (3.6) provides all the information we need for the logarithmic generating function of the linear forms of $L(t,\cdot,\varepsilon)$. If $\{L(t,\cdot,\varepsilon); t \ge 0\}$ were exponentially tight in $\mathcal{L}^p(\mathbb{R}^d)$, then the upper bound for (3.2) would follow from a standard argument. Unfortunately, this is not the case, so we need the following localization procedure to compactify $L(t,\cdot,\varepsilon)$.

Let $m > 0$ be fixed and let $G_m$ be the discrete subgroup of $\mathbb{R}^d$ that consists of vectors whose coordinates are integer multiples of $m$. Let $T_m$ be the quotient of $\mathbb{R}^d$ modulo $G_m$ and let $\iota : \mathbb{R}^d \to T_m$ be the canonical map. Then the $T_m$-valued process $W^*(t) = \iota(W(t))$ is a Markov process and is called, in the literature, Brownian motion on the torus $T_m$. Notice that $T_m$ becomes a compact group under the induced distance

$$d(x^*, y^*) = \inf\{|x-y|; \iota(x) = x^* \text{ and } \iota(y) = y^*\}.$$

Let $\lambda(dx^*)$ be the Lebesgue (Haar) measure on $T_m$ and write

$$\begin{aligned} L^*(t, x^*, \varepsilon) &= \sum_{\mathbf{k}\in\mathbb{Z}^d} L(t, x + m\mathbf{k}, \varepsilon) \\ &= \int_0^t \varphi_\varepsilon(W^*(s) - x^*)\, ds, \qquad t \geq 0, x^* \in T_m, \end{aligned}$$

where $\varphi_\varepsilon$ is a function on $T_m$ defined by

$$\varphi_\varepsilon(x^*) = \frac{1}{C_d \varepsilon^d} \sum_{\mathbf{k}\in\mathbb{Z}^d} \mathbb{1}_{\{|x+m\mathbf{k}|\leq\varepsilon\}}.$$

Notice that if $m$ is large enough, then the above summation has at most one nonzero term, so $\varphi_\varepsilon(x^*) \leq C_d^{-1}\varepsilon^{-d}$ and

$$\begin{aligned} \int_{\mathbb{R}^d} L^p(t, x, \varepsilon)\, dx &= \sum_{\mathbf{k}\in\mathbb{Z}^d} \int_{[0,m]^d} L^p(t, x + m\mathbf{k}, \varepsilon)\, dx \\ &\leq \int_{[0,m]^d} \left[\sum_{\mathbf{k}\in\mathbb{Z}^d} L(t, x + m\mathbf{k}, \varepsilon)\right]^p dx \qquad (3.7) \\ &= \int_{T_m} [L^*(t, x^*, \varepsilon)]^p \lambda(dx^*). \end{aligned}$$

For each integrable function $f$ on $T_m$, define the function $f^*$ on $\mathbb{R}^d$ by $f^* = f \circ \iota$. We have

$$\begin{aligned} \int_{T_m} f(x^*) L^*(t, x^*, \varepsilon) \lambda(dx^*) &= \sum_{\mathbf{k}\in\mathbb{Z}^d} \int_{[0,m]^d} f^*(x) L(t, x + m\mathbf{k}, \varepsilon)\, dx \\ &= \int_{\mathbb{R}^d} f^*(x) L(t, x, \varepsilon)\, dx = \int_0^t (f^*)_\varepsilon(W(s))\, ds. \end{aligned}$$

By (3.6),

$$\begin{aligned} &\lim_{t\to\infty} \frac{1}{t} \log \mathbb{E} \exp\left\{\theta \int_{T_m} f(x^*) L^*(t, x^*, \varepsilon) \lambda(dx^*)\right\} \\ &= \sup_{g\in\mathcal{F}} \left\{\theta \int_{\mathbb{R}^d} (f^*)_\varepsilon(x) g^2(x)\, dx - \frac{1}{2}\int_{\mathbb{R}^d} |\nabla g(x)|^2\, dx\right\} \qquad (3.8) \\ &= \sup_{g\in\mathcal{F}} \left\{\theta \int_{[0,m]^d} f^*(x) \left(\sum_{\mathbf{k}\in\mathbb{Z}^d} (g^2)_\varepsilon(x + m\mathbf{k})\right) dx - \frac{1}{2}\int_{\mathbb{R}^d} |\nabla g(x)|^2\, dx\right\}. \end{aligned}$$

We now intend to establish exponential tightness for $L^*(t,\cdot,\varepsilon)$ as it is viewed as a process taking values in $\mathcal{L}^p(T_m)$. For any $x^*, y^* \in T_m$ with $d(x^*, y^*) \le \delta$, there are $x, y \in \mathbb{R}^d$ such that $\iota(x) = x^*$, $\iota(y) = y^*$ and $|x - y| \le \delta$. Therefore,

$$\begin{aligned}
&|L^*(1, x^*, \varepsilon) - L^*(1, y^*, \varepsilon)| \\
&\quad \le \sum_{\mathbf{k}\in\mathbb{Z}^d} |L(1, x + m\mathbf{k}, \varepsilon) - L(1, y + m\mathbf{k}, \varepsilon)| \\
&\quad \le \#\left\{\mathbf{k}; \min\{|x + m\mathbf{k}|, |y + m\mathbf{k}|\} \le \max_{0\le t\le 1} |W(t)| + \varepsilon\right\} \\
&\qquad \times \sup_{|x-y|\le\delta} |L(1, x, \varepsilon) - L(1, y, \varepsilon)|.
\end{aligned}$$

By shifting invariance, the quantity $\#\{\mathbf{k}\}$ on the right-hand side can be bounded by a finite random number independent of $x$, $y$ and $\delta$. By the obvious fact that

$$\sup_{|x-y|\le\delta} |L(1, x, \varepsilon) - L(1, y, \varepsilon)| \xrightarrow{P} 0, \qquad (\delta \to 0),$$

we have

$$\gamma_\delta \equiv \sup_{d(x^*,y^*)\le\delta} |L^*(1, x^*, \varepsilon) - L^*(1, y^*, \varepsilon)| \xrightarrow{P} 0, \qquad (\delta \to 0).$$

Consequently,

$$\begin{aligned}
&\sup_{d(y^*,z^*)\le\delta} \int_{T_m} |L^*(1, y^* + x^*, \varepsilon) - L^*(1, z^* + x^*, \varepsilon)|^p \lambda(\,dx^*) \\
&\quad \le \gamma_\delta^{p-1} \sup_{d(y^*,z^*)\le\delta} \int_{T_m} \{L^*(1, y^* + x^*, \varepsilon) + L^*(1, z^* + x^*, \varepsilon)\}\lambda(\,dx^*) \\
&\quad = 2\gamma_\delta^{p-1} \xrightarrow{P} 0 \qquad (\delta \to 0).
\end{aligned}$$

Hence, the family $\{L^*(1, y + \cdot, \varepsilon)\}_{y\in T_m}$ of $\mathcal{L}^p(T_m)$-valued random variables is uniformly tight. In addition, for any $y \in T_m$,

$$\begin{aligned}
\int_{T_m} [L^*(1, y^* + x^*, \varepsilon)]^p \lambda(\,dx^*) &\le \left(\frac{1}{C_d\varepsilon^d}\right)^{p-1} \int_{T_m} L^*(1, y^* + x^*, \varepsilon)\lambda(\,dx^*) \\
&= \left(\frac{1}{C_d\varepsilon^d}\right)^{p-1} < \infty.
\end{aligned}$$

By Theorem 3.1 in de Acosta (1985), there is a compact, convex, positively balanced subset $K \subset \mathcal{L}^p(T_m)$ such that

$$\sup_{y\in T_m} \mathbb{E}_y \exp\{q_K(L^*(1, \cdot, \varepsilon))\} < \infty,$$

where $q_K(\cdot)$ is Minkowski functional, which is a seminorm on $\mathcal{L}^p(T_m)$. Applying the Markov property, we have

$$\mathbb{E}\exp\{q_K(L^*(t,\cdot,\varepsilon))\} \le \Big(\sup_{y\in T_m}\mathbb{E}_y\exp\{q_K(L^*(1,\cdot,\varepsilon))\}\Big)^t,$$

which gives

$$\limsup_{t\to\infty}\frac{1}{t}\log\mathbb{E}\exp\{q_K(L^*(t,\cdot,\varepsilon))\}<\infty. \tag{3.9}$$

Notice that for any $\gamma>0$,

$$\begin{aligned}\mathbb{E}\Big[\exp\Big\{\theta\Big(\int_{T_m}[L^*(t,x^*,\varepsilon)]^p\lambda(dx^*)\Big)^{1/p}\Big\};\frac{1}{t}L^*(t,\cdot,\varepsilon)\notin\gamma K\Big]\\ \le \exp\Big\{\theta\Big(\frac{1}{C_d\varepsilon}\Big)^{p-1}t\Big\}\mathbb{P}\{q_K(L^*(t,\cdot,\varepsilon))\ge\gamma t\}.\end{aligned}$$

In view of (3.9), we have that for sufficiently large $\gamma$,

$$\begin{aligned}&\mathbb{E}\exp\Big\{\theta\Big(\int_{T_m}[L^*(t,x^*,\varepsilon)]^p\lambda(dx^*)\Big)^{1/p}\Big\}\\ &\quad\sim\mathbb{E}\Big[\exp\Big\{\theta\Big(\int_{T_m}[L^*(t,x^*,\varepsilon)]^p\lambda(dx^*)\Big)^{1/p}\Big\};\frac{1}{t}L^*(t,\cdot,\varepsilon)\in\gamma K\Big]\end{aligned}\tag{3.10}$$

as $t\to\infty$.

Let $q>1$ be the conjugate number of $p$ and let $\delta>0$ be fixed. By the Hahn–Banach theorem and compactness, there are finitely many bounded functions $f_1,\dots,f_N$ in the unit sphere of $\mathcal{L}^q(T_m)$ such that

$$\Big(\int_{T_m}|h(x^*)|^p\lambda(dx^*)\Big)^{1/p}<\max_{1\le i\le N}\int_{T_m}f_i(x^*)h(x^*)\lambda(dx^*)+\delta\qquad\forall\, h\in\gamma K.$$

In particular,

$$\begin{aligned}\mathbb{E}\Big(\exp\Big\{\theta\Big(\int_{T_m}[L^*(t,x^*,\varepsilon)]^p\lambda(d\bar{x})\Big)^{1/p}\Big\};\frac{1}{t}L^*(t,\cdot,\varepsilon)\in\gamma K\Big)\\ \le e^{\theta\delta t}\sum_{i=1}^N\mathbb{E}\exp\Big\{\theta\int_{T_m}f_i(x^*)L^*(t,x^*,\varepsilon)\lambda(dx^*)\Big\}.\end{aligned}$$

By (3.8) and (3.10),

$$\begin{aligned}&\limsup_{t\to\infty}\frac{1}{t}\log\mathbb{E}\exp\Big\{\theta\Big(\int_{T_m}[L^*(t,x^*,\varepsilon)]^p\lambda(d\bar{x})\Big)^{1/p}\Big\}\\ &\quad\le\theta\delta+\max_{1\le i\le N}\sup_{g\in\mathcal{F}_d}\Big\{\theta\int_{[0,m]^d}f_i^*(x)\Big(\sum_{\mathbf{k}\in\mathbb{Z}^d}(g^2)_\varepsilon(x+m\mathbf{k})\Big)dx\end{aligned}$$

$$
\begin{aligned}
&\qquad\qquad\qquad\qquad - \frac{1}{2}\int_{\mathbb{R}^d} |\nabla g(x)|^2\,dx\bigg\} \\
&\le \theta\delta + \sup_{g\in\mathcal{F}_d}\bigg\{\theta\bigg(\int_{[0,m]^d}\bigg(\sum_{\mathbf{k}\in\mathbb{Z}^d}(g^2)_\varepsilon(x+m\mathbf{k})\bigg)^p dx\bigg)^{1/p} \\
&\qquad\qquad\qquad\qquad - \frac{1}{2}\int_{\mathbb{R}^d} |\nabla g(x)|^2\,dx\bigg\},
\end{aligned}
$$

where the last step follows from Hölder's inequality. In view of (3.7) and Lemma A.1, letting $\delta \to 0$ and then $m \to \infty$, we obtain the upper bound for (3.2):

$$
\limsup_{t\to\infty} \frac{1}{t}\log \mathbb{E}\exp\bigg\{\theta\bigg(\int_{\mathbb{R}^d} L^p(t,x,\varepsilon)\,dx\bigg)^{1/p}\bigg\} \le M_\varepsilon(\theta). \tag{3.11}
$$

By the inequality

$$
\bigg(\int_{\mathbb{R}^d}\prod_{j=1}^p L_j(t,x,\varepsilon)\,dx\bigg)^{1/p} \le \frac{1}{p}\sum_{j=1}^p\bigg(\int_{\mathbb{R}^d} L_j^p(t,x,\varepsilon)\,dx\bigg)^{1/p}
$$

we have

$$
\mathbb{E}\exp\bigg\{\theta\bigg(\int_{\mathbb{R}^d}\prod_{j=1}^p L_j(t,x,\varepsilon)\,dx\bigg)^{1/p}\bigg\} \le \bigg[\mathbb{E}\exp\bigg\{\frac{\theta}{p}\bigg(\int_{\mathbb{R}^d} L^p(t,x,\varepsilon)\,dx\bigg)^{1/p}\bigg\}\bigg]^p.
$$

From (3.11) (with $\theta$ replaced by $p^{-1}\theta$) we have the upper bound for (3.3):

$$
\limsup_{t\to\infty} \frac{1}{t}\log \mathbb{E}\exp\bigg\{\theta\bigg(\int_{\mathbb{R}^d}\prod_{j=1}^p L_j(t,x,\varepsilon)\,dx\bigg)^{1/p}\bigg\} \le N_\varepsilon(\theta). \tag{3.12}
$$

We now come to the proof of the lower bounds for (3.2) and (3.3). Notice that for any $0 < r \le \infty$ and any measurable function $f$ on $\mathbb{R}^d$ with $\|f\|_q = 1$ and $f(x) = 0$ for $|x| > r$ (if $r < \infty$),

$$
\bigg(\int_{\{|x|\le r\}} L^p(t,x,\varepsilon)\,dx\bigg)^{1/p} \ge \int_{\mathbb{R}^d} f(x)L(t,x,\varepsilon)\,dx = \int_0^t f_\varepsilon(W(s))\,ds.
$$

By (3.6) we have

$$
\begin{aligned}
&\liminf_{t\to\infty}\frac{1}{t}\log\mathbb{E}\exp\bigg\{\theta\bigg(\int_{\{|x|\le r\}} L^p(t,x,\varepsilon)\,dx\bigg)^{1/p}\bigg\} \\
&\quad\ge \sup_{g\in\mathcal{F}_d}\bigg\{\theta\int_{\mathbb{R}^d} f_\varepsilon(x)g^2(x)\,dx - \frac{1}{2}\int_{\mathbb{R}^d}|\nabla g(x)|^2\,dx\bigg\} \\
&\quad= \sup_{g\in\mathcal{F}_d}\bigg\{\theta\int_{\{|x|\le r\}} f(x)(g^2)_\varepsilon(x)\,dx - \frac{1}{2}\int_{\mathbb{R}^d}|\nabla g(x)|^2\,dx\bigg\}.
\end{aligned}
$$

Taking the supremum over $f$ on the right-hand side,

$$\begin{aligned}(3.13)\qquad &\liminf_{t\to\infty}\frac{1}{t}\log\mathbb{E}\exp\left\{\theta\left(\int_{\{|x|\le r\}}L^p(t,x,\varepsilon)\,dx\right)^{1/p}\right\}\\ &\quad\ge\sup_{g\in\mathcal{F}_d}\left\{\theta\left(\int_{\{|x|\le r\}}|(g^2)_\varepsilon(x)|^p\,dx\right)^{1/p}-\frac{1}{2}\int_{\mathbb{R}^d}|\nabla g(x)|^2\,dx\right\}.\end{aligned}$$

In particular, letting $r=\infty$ gives the lower bound for (3.2).

To prove the lower bound for (3.3), we view $L(t,\cdot,\varepsilon)$ as a process with values in $\mathcal{L}^p(B_r)$ by limiting $x$ to $B_r$, where $B_r$ is the $d$-dimensional closed ball with center 0 and radius $r>0$ being fixed but arbitrary. We need to show that $t^{-1}L(t,\cdot,\varepsilon)$ is exponentially tight in $\mathcal{L}^p(B_r)$: For any $\gamma>0$, there is a compact set $K_0$ in $\mathcal{L}^p(B_r)$, such that

$$(3.14)\qquad \limsup_{t\to\infty}\frac{1}{t}\log\mathbb{P}\left\{\frac{1}{t}L(t,\cdot,\varepsilon)\notin K_0\right\}\le-\gamma.$$

To this end, we first prove that for any $\lambda>0$,

$$\begin{aligned}(3.15)\qquad &\lim_{\delta\to0^+}\sup_{t\ge1}\frac{1}{t}\log\mathbb{E}\exp\left\{\lambda\sup_{|y|\le\delta}\left(\int_{\mathbb{R}^d}|L(t,x+y,\varepsilon)-L(t,x,\varepsilon)|^p\,dx\right)^{1/p}\right\}\\ &\quad=0.\end{aligned}$$

Indeed, by subadditivity,

$$\begin{aligned}&\mathbb{E}\exp\left\{\lambda\sup_{|y|\le\delta}\left(\int_{\mathbb{R}^d}|L(t,x+y,\varepsilon)-L(t,x,\varepsilon)|^p\,dx\right)^{1/p}\right\}\\ &\quad\le\left[\mathbb{E}\exp\left\{\lambda\sup_{|y|\le\delta}\left(\int_{\mathbb{R}^d}|L(1,x+y,\varepsilon)-L(1,x,\varepsilon)|^p\,dx\right)^{1/p}\right\}\right]^t.\end{aligned}$$

To have (3.15), we need only to show

$$(3.16)\qquad \lim_{\delta\to0^+}\mathbb{E}\exp\left\{\lambda\sup_{|y|\le\delta}\left(\int_{\mathbb{R}^d}|L(1,x+y,\varepsilon)-L(1,x,\varepsilon)|^p\,dx\right)^{1/p}\right\}=1.$$

Notice that as $|y|\le\delta<\varepsilon$,

$$\begin{aligned}&|L(1,x+y,\varepsilon)-L(1,x,\varepsilon)|\\ &\quad\le\frac{1}{C_d\varepsilon^d}\left\{\int_0^1\mathbb{1}_{\{\varepsilon-\delta<|W(s)-x|\le\varepsilon\}}\,ds+\int_0^1\mathbb{1}_{\{\varepsilon-\delta<|W(s)-(x+y)|\le\varepsilon\}}\,ds\right\}\\ &\quad=\xi_\delta(x)+\xi_\delta(x+y)\qquad\text{(say)}.\end{aligned}$$

By triangular inequality,

$$\sup_{|y|\le\delta}\left(\int_{\mathbb{R}^d}|L(1,x+y,\varepsilon)-L(1,x,\varepsilon)|^p\,dx\right)^{1/p}\le2\left(\int_{\mathbb{R}^d}|\xi_\delta(x)|^p\,dx\right)^{1/p}.$$

The right-hand side approaches 0 as $\delta \to 0^+$. Hence (3.16) follows from the dominated convergence theorem.

For each $k \geq 1$, by (3.15) and Chebyshev's inequality there is a $\delta_k > 0$ such that

$$\mathbb{P}\left\{\sup_{|y|\leq\delta_k}\left(\int_{\mathbb{R}^d}|L(t,x+y,\varepsilon)-L(t,x,\varepsilon)|^p\,dx\right)^{1/p}\geq k^{-1}t\right\}\leq e^{-k\gamma t}\qquad\forall\, t\geq 1.$$

In addition, by (3.12) there is a $C > 0$ such that

$$\limsup_{t\to\infty}\frac{1}{t}\log\mathbb{P}\left\{\left(\int_{\mathbb{R}^d}L^p(t,x,\varepsilon)\,dx\right)^{1/p}\geq Ct\right\}\leq -\gamma.$$

Consider the set $A \subset \mathcal{L}^p(B_r)$ given by

$$A=\left\{f\in\mathcal{L}^p(B_r);\|f\|_p\leq C \text{ and } \sup_{|y|\leq\delta_k}\|f(\cdot+y)-f(\cdot)\|_p\leq k^{-1}\ \forall\, k\geq 1\right\}.$$

By the criterion of compactness in $L_p$ space [see, e.g., Dunford and Schwartz (1988), Theorem 21, page 301], $A$ is conditionally compact in $\mathcal{L}^p(B_r)$. In addition, from the construction of $A$ we have

$$\limsup_{t\to\infty}\frac{1}{t}\log\mathbb{P}\left\{\frac{1}{t}L(t,\cdot,\varepsilon)\notin A\right\}\leq -\gamma.$$

Taking $K_0$ as the closure of $A$, we have (3.14).

Taking $\gamma$ sufficiently large, we have

$$\begin{aligned}&\mathbb{E}\exp\left\{\frac{\theta}{p}\left(\int_{\{|x|\leq r\}}L^p(t,x,\varepsilon)\,dx\right)^{1/p}\right\}\\ &\quad\sim\mathbb{E}\left[\exp\left\{\frac{\theta}{p}\left(\int_{\{|x|\leq r\}}L^p(t,x,\varepsilon)\,dx\right)^{1/p}\right\};\frac{1}{t}L(t,\cdot,\varepsilon)\in K_0\right].\end{aligned}\tag{3.17}$$

Consider the continuous, nonnegative functional $\Psi$ defined on $(\mathcal{L}^p(B_r))^p$:

$$\Psi(f_1,\ldots,f_p)=\frac{1}{p}\sum_{j=1}^p\left(\int_{\{|x|\leq r\}}|f_j(x)|^p\,dx\right)^{1/p}-\left(\int_{\{|x|\leq r\}}\prod_{j=1}^p|f_j(x)|\,dx\right)^{1/p}.$$

Clearly, $\Psi \equiv 0$ on the diagonal

$$\{(f_1,\ldots,f_p);f_1=\cdots=f_p\}.$$

Hence, for given $\delta > 0$ and any $g \in \mathcal{L}^p(B_r)$, there exists a $b = b(g,\delta) > 0$ such that

$$\Psi(f_1,\ldots,f_p)\leq\delta\qquad\text{as } f_j\in B(g,b)\ \forall\, 1\leq j\leq p,$$

where $B(g,b)$ stands for the open ball in $\mathcal{L}^p(B_r)$ with center $g$ and radius $b$. Therefore,

$$\mathbb{E}\exp\biggl\{\theta\biggl(\int_{\{|x|\le r\}}\prod_{j=1}^p L_j(t,x,\varepsilon)\,dx\biggr)^{1/p}\biggr\}$$

$$\ge e^{-\delta t}\mathbb{E}\biggl[\exp\biggl\{\frac{\theta}{p}\sum_{j=1}^p\biggl(\int_{\{|x|\le r\}}L_j^p(t,x,\varepsilon)\,dx\biggr)^{1/p}\biggr\};\ \frac{1}{t}L_j(t,\cdot,\varepsilon)\in B(g,b)\ \forall\, 1\le j\le p\biggr] \tag{3.18}$$

$$= e^{-\delta t}\biggl(\mathbb{E}\biggl[\exp\biggl\{\frac{\theta}{p}\biggl(\int_{\{|x|\le r\}}L^p(t,x,\varepsilon)\,dx\biggr)^{1/p}\biggr\};\frac{1}{t}L(t,\cdot,\varepsilon)\in B(g,b)\biggr]\biggr)^p.$$

Let $\{B(g_1,b_1),\dots,B(g_N,b_N)\}$ be a finite subfamily of the open sets

$$\{B(g,b(g,\delta));\ g\in K_0\},$$

which cover $K_0$. Then

$$\mathbb{E}\biggl[\exp\biggl\{\frac{\theta}{p}\biggl(\int_{\{|x|\le r\}}L^p(t,x,\varepsilon)\,dx\biggr)^{1/p}\biggr\};\frac{1}{t}L(t,\cdot,\varepsilon)\in K_0\biggr]$$

$$\le\sum_{i=1}^N\mathbb{E}\biggl[\exp\biggl\{\frac{\theta}{p}\biggl(\int_{B_r}L^p(t,x,\varepsilon)\,dx\biggr)^{1/p}\biggr\};\frac{1}{t}L(t,\cdot,\varepsilon)\in B(g_i,b_i)\biggr].$$

Therefore,

$$\liminf_{t\to\infty}\frac{1}{t}\log\max_{1\le i\le N}\mathbb{E}\biggl[\exp\biggl\{\frac{\theta}{p}\biggl(\int_{\{|x|\le r\}}L^p(t,x,\varepsilon)\,dx\biggr)^{1/p}\biggr\};\ \frac{1}{t}L(t,\cdot,\varepsilon)\in B(g_i,b_i)\biggr]$$

$$\ge\liminf_{t\to\infty}\frac{1}{t}\log\mathbb{E}\biggl[\exp\biggl\{\frac{\theta}{p}\biggl(\int_{\{|x|\le r\}}L^p(t,x,\varepsilon)\,dx\biggr)^{1/p}\biggr\};\frac{1}{t}L(t,\cdot,\varepsilon)\in K_0\biggr].$$

Combining this with (3.13) (with $\theta$ replaced by $p^{-1}\theta$), (3.17) and (3.18), we have

$$\liminf_{t\to\infty}\frac{1}{t}\log\mathbb{E}\exp\biggl\{\theta\biggl(\int_{\{|x|\le r\}}\prod_{j=1}^p L_j(t,x,\varepsilon)\,dx\biggr)^{1/p}\biggr\}$$

$$\ge -\delta + p\sup_{g\in\mathcal{F}}\biggl\{\frac{\theta}{p}\biggl(\int_{\{|x|\le r\}}|(g^2)_\varepsilon(x)|^p\,dx\biggr)^{1/p}-\frac{1}{2}\int_{\mathbb{R}^d}|\nabla g(x)|^2\,dx\biggr\}.$$

Letting $\delta \to 0^+$ and $r \to \infty$, we obtain the lower bound for (3.3):

$$\liminf_{t\to\infty} \frac{1}{t} \log \mathbb{E} \exp\left\{\theta\left(\int_{\mathbb{R}^d} \prod_{j=1}^p L_j(t,x,\varepsilon)\, dx\right)^{1/p}\right\} \ge N_\varepsilon(\theta). \tag{3.19}$$

We now come to the random walks. Given $t > 0$, write $t_n = [tn/b_n]$ and $\gamma_n = [n/t_n]$. Then $n \le t_n(\gamma_n + 1)$. By independence and triangular inequality,

$$\begin{aligned}
&\mathbb{E}\exp\left\{\theta\left(\frac{b_n}{n}\right)^{(2p-d(p-1))/(2p)}\left(\sum_{x\in\mathbb{Z}^d} l^p(n,x,\varepsilon)\right)^{1/p}\right\} \\
&\quad\le \Bigg(\mathbb{E}\exp\Bigg\{\theta\left(\frac{b_n}{n}\right)^{(2p-d(p-1))/(2p)} \\
&\qquad\qquad \times\left(\sum_{x\in\mathbb{Z}^d}\left[\frac{1}{\#\{B_n\}}\sum_{k=1}^{t_n}\mathbb{1}\left\{\left|\frac{S(k)-x}{\sqrt{nb_n^{-1}}}\right|\le\varepsilon\right\}\right]^p\right)^{1/p}\Bigg\}\Bigg)^{\gamma_n+1} \\
&\quad= \Bigg(\mathbb{E}\exp\Bigg\{\theta\left(\frac{b_n}{n}\right)^{(2p-d(p-1))/(2p)} \\
&\qquad\qquad \times\left(\int_{\mathbb{R}^d}\left[\frac{1}{\#\{B_n\}}\sum_{k=1}^{t_n}\mathbb{1}\left\{\left|\frac{S(k)-[x]}{\sqrt{nb_n^{-1}}}\right|\le\varepsilon\right\}\right]^p dx\right)^{1/p}\Bigg\}\Bigg)^{\gamma_n+1} \\
&\quad= \Bigg(\mathbb{E}\exp\Bigg\{\theta\Bigg(\int_{\mathbb{R}^d}\Bigg[\left(\frac{b_n}{n}\right)^{(2-d)/2} \\
&\qquad\qquad \times\frac{1}{\#\{B_n\}}\sum_{k=1}^{t_n}\mathbb{1}\left\{\left|\frac{S(k)-[\sqrt{nb_n^{-1}}x]}{\sqrt{nb_n^{-1}}}\right|\le\varepsilon\right\}\Bigg]^p dx\Bigg)^{1/p}\Bigg\}\Bigg)^{\gamma_n+1}.
\end{aligned}$$

Notice that

$$\left(\frac{b_n}{n}\right)^{(2-d)/2}\frac{1}{\#\{B_n\}} \sim \frac{b_n}{n} C_d^{-1}\varepsilon^{-d} \qquad (n\to\infty).$$

Applying the invariance principle to the continuous, uniformly bounded and uniformly convergent functionals $\{\varphi_n(f)\}$ on $C\{[0,t];\mathbb{R}^d\}$, given as

$$\varphi_n(f) = \int_{\mathbb{R}^d}\left(\frac{1}{C_d\varepsilon^d}\int_0^t \mathbb{1}\left\{\left|f(s) - \frac{[\sqrt{nb_n^{-1}}x]}{\sqrt{nb_n^{-1}}}\right|\le\varepsilon\right\}ds\right)^p dx,$$

we have

$$\lim_{n\to\infty}\mathbb{E}\exp\left\{\theta\left(\int_{\mathbb{R}^d}\left[\left(\frac{b_n}{n}\right)^{(2-d)/2}\right.\right.\right.$$

$$\times \frac{1}{\#\{B_n\}} \sum_{k=1}^{t_n} \mathbb{1}\left\{\left|\frac{S(k)-[\sqrt{nb_n^{-1}}x]}{\sqrt{nb_n^{-1}}}\right| \le \varepsilon\right\}\Bigg]^p dx\Bigg)^{1/p}\Bigg\}$$

$$= \mathbb{E}\exp\left\{\theta\left(\int_{\mathbb{R}^d} \widetilde{L}^p(t,x,\varepsilon)\,dx\right)^{1/p}\right\}$$

where $\widetilde{L}(t,x,\varepsilon)$ is the analogue of $L(t,x,\varepsilon)$ with $W$ replaced by the Gaussian Lévy process $U$ whose covariance matrix is $\Gamma$. Therefore,

$$\begin{aligned} &\limsup_{n\to\infty} \frac{1}{b_n} \log \mathbb{E}\exp\left\{\theta\left(\frac{b_n}{n}\right)^{(2p-d(p-1))/(2p)} \left(\sum_{x\in\mathbb{Z}^d} l^p(n,x,\varepsilon)\right)^{1/p}\right\} \\ &\qquad \le \frac{1}{t}\log \mathbb{E}\exp\left\{\theta\left(\int_{\mathbb{R}^d} \widetilde{L}^p(t,x,\varepsilon)\,dx\right)^{1/p}\right\}. \end{aligned} \tag{3.20}$$

By the same argument used in the canonical Brownian case,

$$\lim_{t\to\infty} \frac{1}{t}\log \mathbb{E}\exp\left\{\theta\left(\int_{\mathbb{R}^d} \widetilde{L}^p(t,x,\varepsilon)\,dx\right)^{1/p}\right\} = \widetilde{M}_\varepsilon(\theta).$$

Hence, letting $t\to\infty$ gives the upper bound for (3.4):

$$\begin{aligned} &\limsup_{n\to\infty} \frac{1}{b_n} \log \mathbb{E}\exp\left\{\theta\left(\frac{b_n}{n}\right)^{(2p-d(p-1))/(2p)} \left(\sum_{x\in\mathbb{Z}^d} l^p(n,x,\varepsilon)\right)^{1/p}\right\} \\ &\qquad \le \widetilde{M}_\varepsilon(\theta). \end{aligned} \tag{3.21}$$

By the inequality

$$\left(\sum_{x\in\mathbb{Z}^d}\prod_{j=1}^p l_j(n,x,\varepsilon)\right)^{1/p} \le \frac{1}{p}\sum_{j=1}^p \left(\sum_{x\in\mathbb{Z}^d} l_j^p(n,x,\varepsilon)\right)^{1/p}$$

and (3.21), we have the upper bound for (3.5):

$$\begin{aligned} &\limsup_{n\to\infty} \frac{1}{b_n} \log \mathbb{E}\exp\left\{\theta\left(\frac{b_n}{n}\right)^{(2p-d(p-1))/(2p)} \left(\sum_{x\in\mathbb{Z}^d}\prod_{j=1}^p l_j(n,x,\varepsilon)\right)^{1/p}\right\} \\ &\qquad \le \widetilde{M}_\varepsilon(\theta). \end{aligned} \tag{3.22}$$

On the other hand, for any uniformly continuous function $f$ supported on $B_r$ with $\|f\|_q = 1$,

$$\begin{aligned} \left(\sum_{x\in\mathbb{Z}^d} l^p(n,x,\varepsilon)\,dx\right)^{1/p} &= \left(\frac{n}{b_n}\right)^{d/(2p)} \left(\int_{\mathbb{R}^d} l^p(n,[\sqrt{nb_n^{-1}}x],\varepsilon)\,dx\right)^{1/p} \\ &\ge \left(\frac{n}{b_n}\right)^{d/(2p)} \int_{\mathbb{R}^d} f(x)l(n,[\sqrt{nb_n^{-1}}x],\varepsilon)\,dx. \end{aligned}$$

Notice that

$$
\begin{aligned}
&\int_{\mathbb{R}^d} f(x) l(n, [\sqrt{n b_n^{-1}} x], \varepsilon)\, dx \\
&\quad = \left(\frac{b_n}{n}\right)^{d/2} \int_{\mathbb{R}^d} f\left(\sqrt{\frac{b_n}{n}} x\right) l(n, [x], \varepsilon)\, dx \\
&\quad = \left(\frac{b_n}{n}\right)^{d/2} \left\{ o(n) + \sum_{x \in \mathbb{Z}^d} f\left(\sqrt{\frac{b_n}{n}} x\right) l(n, x, \varepsilon) \right\} \\
&\quad = \left(\frac{b_n}{n}\right)^{d/2} \left\{ o(n) + \sum_{k=1}^{n} \frac{1}{\#\{B_n\}} \sum_{x \in \mathbb{Z}^d} f\left(\sqrt{\frac{b_n}{n}} (x + S(k))\right) \mathbb{1}_{\{x \in B_n\}} \right\} \\
&\quad = \left(\frac{b_n}{n}\right)^{d/2} \left\{ o(n) + \sum_{k=1}^{n} f_\varepsilon\left(\sqrt{\frac{b_n}{n}} S(k)\right) \right\} \qquad (n \to \infty),
\end{aligned}
$$

where the term $o(n)$ is bounded by a deterministic sequence $a_n$ that satisfies $a_n/n \to 0$ as $n \to \infty$. By Theorem 4.1 in Chen and Li (2004),

$$
\begin{aligned}
&\liminf_{n\to\infty} \frac{1}{b_n} \log \mathbb{E} \exp\left\{ \theta \left(\frac{b_n}{n}\right)^{(2-d)/2} \int_{\mathbb{R}^d} f(x) l(n, [\sqrt{n b_n^{-1}} x], \varepsilon)\, dx \right\} \\
&\qquad \geq \sup_{g \in \mathcal{F}_d} \left\{ \theta \int_{\mathbb{R}^d} f_\varepsilon(x) g^2(x)\, dx - \frac{1}{2} \int_{\mathbb{R}^d} \langle \nabla g(x), \Gamma \nabla g(x) \rangle\, dx \right\} \\
&\qquad = \sup_{g \in \mathcal{F}_d} \left\{ \theta \int_{\mathbb{R}^d} f(x) (g^2)_\varepsilon(x)\, dx - \frac{1}{2} \int_{\mathbb{R}^d} \langle \nabla g(x), \Gamma \nabla g(x) \rangle\, dx \right\}.
\end{aligned} \tag{3.23}
$$

Consequently, we have the lower bound for (3.4):

$$
\begin{aligned}
&\liminf_{n\to\infty} \frac{1}{b_n} \log \mathbb{E} \exp\left\{ \theta \left(\frac{b_n}{n}\right)^{(2p-d(p-1))/(2p)} \left( \sum_{x \in \mathbb{Z}^d} l^p(n, x, \varepsilon) \right)^{1/p} \right\} \\
&\qquad \geq \widetilde{M}_\varepsilon(\theta).
\end{aligned} \tag{3.24}
$$

It remains to prove the lower bound for (3.5). Notice that

$$
\sum_{x \in \mathbb{Z}^d} \prod_{j=1}^{p} l_j(n, x, \varepsilon) = \left(\frac{n}{b_n}\right)^{d/2} \int_{\mathbb{R}^d} \prod_{j=1}^{p} l_j(n, [\sqrt{n b_n^{-1}} x], \varepsilon)\, dx.
$$

We need only to prove

$$
\begin{aligned}
&\liminf_{n\to\infty} \frac{1}{b_n} \log \mathbb{E} \exp\left\{ \theta \left(\frac{b_n}{n}\right)^{(2-d)/2} \left( \int_{\mathbb{R}^d} \prod_{j=1}^{p} l_j(n, [\sqrt{n b_n^{-1}} x], \varepsilon)\, dx \right)^{1/p} \right\} \\
&\qquad \geq \widetilde{N}_\varepsilon(\theta).
\end{aligned} \tag{3.25}
$$

Similar to (3.20) (with $t=1$), we can prove that for any $\theta>0$ and $\delta>0$,

$$\begin{aligned}\limsup_{n\to\infty}\frac{1}{b_n}\log\mathbb{E}\exp\biggl\{\theta\biggl(\frac{b_n}{n}\biggr)^{(2-d)/2}&\\ \times\sup_{|y|\le\delta}\biggl(\int_{\mathbb{R}^d}|l(n,[\sqrt{nb_n^{-1}}(x+y)],\varepsilon)&\\ -l(n,[\sqrt{nb_n^{-1}}x],\varepsilon)|^p\,dx\biggr)^{1/p}\biggr\}&\\ \le\log\mathbb{E}\exp\biggl\{\theta\sup_{|y|\le\delta}\biggl(\int_{\mathbb{R}^d}|\widetilde{L}(1,x+y,\varepsilon)-\widetilde{L}(1,x,\varepsilon)|^p\,dx\biggr)^{1/p}\biggr\}.&\end{aligned}$$

Replacing $L(t,x,\varepsilon)$ with $\widetilde{L}(t,x,\varepsilon)$ in (3.16), we have

$$\begin{aligned}\lim_{\delta\to0^+}\limsup_{n\to\infty}\frac{1}{b_n}\log\mathbb{E}\exp\biggl\{\theta\biggl(\frac{b_n}{n}\biggr)^{(2-d)/2}&\\ \times\sup_{|y|\le\delta}\biggl(\int_{\mathbb{R}^d}|l(n,[\sqrt{nb_n^{-1}}(x+y)],\varepsilon)&\\ -l(n,[\sqrt{nb_n^{-1}}x],\varepsilon)|^p\,dx\biggr)^{1/p}\biggr\}=0.&\end{aligned}$$

Similar to (3.17), as $x$ is limited to a finite ball $B_r$ and $l(n,[\sqrt{nb_n^{-1}}(\cdot)],\varepsilon)$ is viewed as a process with values in $\mathcal{L}^p(B_r)$, there is a compact set $K_0\subset\mathcal{L}^p(B_r)$ such that

$$\begin{aligned}&\mathbb{E}\exp\biggl\{\frac{\theta}{p}\biggl(\frac{b_n}{n}\biggr)^{(2-d)/2}\biggl(\int_{\mathbb{R}^d}l^p(n,[\sqrt{nb_n^{-1}}x],\varepsilon)\,dx\biggr)^{1/p}\biggr\}\\ &\quad\sim\mathbb{E}\biggl[\exp\biggl\{\frac{\theta}{p}\biggl(\frac{b_n}{n}\biggr)^{(2-d)/2}\biggl(\int_{\mathbb{R}^d}l^p(n,[\sqrt{nb_n^{-1}}x],\varepsilon)\,dx\biggr)^{1/p}\biggr\};\\ &\qquad\qquad n^{(d-2)/2}b_n^{-d/2}l(n,[\sqrt{nb_n^{-1}}(\cdot)],\varepsilon)\in K_0\biggr]\end{aligned}\tag{3.26}$$

as $n\to\infty$. In view of (3.23) and (3.26), the argument used in the proof of (3.16) gives (3.25). $\square$

**4. Upper bound in Theorem 2.1.** In this section, we establish the upper bound for Theorem 2.1:

$$\limsup_{t\to\infty}t^{-1}\log\mathbb{P}\{\alpha([0,1]^p)\ge t^{(d(p-1))/2}\}\le-\frac{p}{2}\kappa(p,d)^{-4p/(d(p-1))}.\tag{4.1}$$

Let $\tau$ be exponential time with parameter 1 and with independent copies $\tau_1, \dots, \tau_p$. Recall that $\mathcal{F}_d$ is defined in (3.1). Write

$$M = \sup_{f \in \mathcal{F}_d} \left\{ \left( \int_{\mathbb{R}^d} |f(x)|^{2p} \, dx \right)^{1/p} - \frac{1}{2} \int_{\mathbb{R}^d} |\nabla f(x)|^2 \, dx \right\}. \tag{4.2}$$

See Lemma 8.2 for how $M$ is related to the Gagliardo–Nirenberg constant.

To apply Theorem 3.1, first notice that by Jensen's inequality,

$$\begin{aligned} M_\varepsilon(\theta) &\le \sup_{f \in \mathcal{F}_d} \left\{ \theta \left( \int_{\mathbb{R}^d} |f(x)|^{2p} \, dx \right)^{1/p} - \frac{1}{2} \int_{\mathbb{R}^d} |\nabla f(x)|^2 \, dx \right\} \\ &= \theta^{2p/(2p-d(p-1))} M, \end{aligned}$$

where the last step follows from the substitution

$$f \mapsto \theta^{dp/(2(2p-d(p-1)))} f(\theta^{p/(2(2p-d(p-1)))} x).$$

By (3.2) in Theorem 3.1,

$$\mathbb{E} \exp\left\{ \theta \left( \int_{\mathbb{R}^d} L^p(\tau, x, \varepsilon) \, dx \right)^{1/p} \right\} < \infty \qquad \forall \, \theta < M^{-(2p-d(p-1))/(2p)}.$$

From the fact

$$\left( \int_{\mathbb{R}^d} \prod_{j=1}^p L_j(\tau_j, x, \varepsilon) \, dx \right)^{1/p} \le \frac{1}{p} \sum_{k=1}^p \left( \int_{\mathbb{R}^d} L_j^p(\tau_j, x, \varepsilon) \, dx \right)^{1/p}$$

we have

$$\begin{aligned} A_\theta &\equiv \mathbb{E} \exp\left\{ \theta \left( \int_{\mathbb{R}^d} \prod_{j=1}^p L_j(\tau_j, x, \varepsilon) \, dx \right)^{1/p} \right\} \\ &\le \left[ \mathbb{E} \exp\left\{ \frac{\theta}{p} \left( \int_{\mathbb{R}^d} L^p(\tau, x, \varepsilon) \, dx \right)^{1/p} \right\} \right]^p \\ &< \infty, \qquad \theta < p \cdot M^{-(2p-d(p-1))/(2p)} \end{aligned}$$

and, therefore,

$$\mathbb{E} \left( \int_{\mathbb{R}^d} \prod_{j=1}^p L_j(\tau_j, x, \varepsilon) \, dx \right)^m \le (pm)! \theta^{-pm} A_\theta \qquad \forall \, \theta < p \cdot M^{-(2p-d(p-1))/(2p)} \tag{4.3}$$

for any integer $m \ge 1$.

On the other hand, write $\delta_x^\varepsilon(y) = (1/C_d\varepsilon^d)\mathbb{1}_{\{|y-x|\le\varepsilon\}}$. For any integer $m \ge 1$, let $\Sigma_m$ be the set of the permutations of $\{1, 2, \ldots, m\}$. Under the convention $\sigma(0) = 0$ and $s_0 = 0$,

$$\begin{aligned}
\mathbb{E}\biggl(\int_{\mathbb{R}^d}\prod_{j=1}^p L_j(\tau_j, x, \varepsilon)\,dx\biggr)^m
&= \int_{(\mathbb{R}^d)^m} dx_1\cdots dx_m\biggl(\mathbb{E}\prod_{k=1}^m L(\tau, x_k, \varepsilon)\biggr)^p \\
&= \int_{(\mathbb{R}^d)^m} dx_1\cdots dx_m\biggl(\sum_{\sigma\in\Sigma_m}\int_0^\infty dt\, e^{-t} \\
&\qquad\times \int_{0\le s_1\le\cdots\le s_m\le t}\mathbb{E}\prod_{k=1}^m \delta_{x_{\sigma(k)}}^\varepsilon(W(s_k))\,ds_1\cdots ds_m\biggr)^p \\
&\ge \int_{(\mathbb{R}^d)^m} dx_1\cdots dx_m\biggl(\sum_{\sigma\in\Sigma_m}\int_0^\infty dt\, e^{-t}\int_{0\le s_1\le\cdots\le s_m\le t} ds_1\cdots ds_m \\
&\qquad\times \prod_{k=1}^m \inf_{|y|\le 2\varepsilon} p_{s_k-s_{k-1}}(x_{\sigma(k)} - x_{\sigma(k-1)} + y)\biggr)^p,
\end{aligned}$$

where the last step follows from the Markov property. Notice that for any $x, y \in \mathbb{R}^d$ with $x \ne 0$ and $|y| \le 2\varepsilon$, and for any $t > 0$,

$$p_t(x+y) \ge p_t\biggl(\biggl(1 + \frac{2\varepsilon}{|x|}\biggr)x\biggr) \equiv q_t(x) \qquad \text{(say).}$$

By (1.9),

$$\begin{aligned}
\mathbb{E}\biggl(\int_{\mathbb{R}^d}\prod_{j=1}^p L_j(\tau_j, x, \varepsilon)\,dx\biggr)^m
&\ge \int_{(\mathbb{R}^d)^m} dx_1\cdots dx_m\biggl(\sum_{\sigma\in\Sigma_m}\int_0^\infty dt\, e^{-t}\int_{0\le s_1\le\cdots\le s_n\le t} ds_1\cdots ds_m \\
&\qquad\times \prod_{k=1}^m q_{s_k-s_{k-1}}(x_{\sigma(k)} - x_{\sigma(k-1)})\biggr)^p \\
&= \int_{(\mathbb{R}^d)^m} dx_1\cdots dx_m\biggl(\sum_{\sigma\in\Sigma_m}\prod_{k=1}^m\int_0^\infty e^{-t} q_t(x_{\sigma(k)} - x_{\sigma(k-1)})\,dt\biggr)^p.
\end{aligned}$$

Let $\delta > 0$ be fixed. We can see that for any $0 < \lambda < 1$, if we take $\varepsilon > 0$ small enough,

$$q_t(x) \ge \lambda p_t(\lambda^{-1}x)$$

holds uniformly for $t \ge \delta$, $x \in \mathbb{R}^d$. Hence,

$$
\begin{aligned}
&\mathbb{E}\left(\int_{\mathbb{R}^d} \prod_{j=1}^p L_j(\tau_j, x, \varepsilon)\, dx\right)^m \\
&\quad \ge \lambda^{pm} \int_{(\mathbb{R}^d)^m} dx_1 \cdots dx_m \left(\sum_{\sigma \in \Sigma_m} \prod_{k=1}^m \int_\delta^\infty e^{-t} p_t(\lambda^{-1}(x_{\sigma(k)} - x_{\sigma(k-1)}))\, dt\right)^p \\
&\quad \ge \lambda^{(p+d)m} e^{-p\delta m} \int_{(\mathbb{R}^d)^m} dx_1 \cdots dx_m \\
&\qquad \times \left(\sum_{\sigma \in \Sigma_m} \prod_{k=1}^m \int_0^\infty e^{-t} p_{t+\delta}(x_{\sigma(k)} - x_{\sigma(k-1)})\, dt\right)^p \\
&\quad = \lambda^{(p+d)m} e^{-p\delta m} \mathbb{E}[\alpha([\delta, \delta + \tau_1] \times \cdots \times [\delta, \delta + \tau_p])^m],
\end{aligned}
$$

where the last step can be seen from the fact (1.9) and Le Gall's [(1990), (1), page 182] moment formula: For any $t_1, \dots, t_p \ge 0$,

$$
\begin{aligned}
&\mathbb{E}\alpha([\delta, \delta + t_1] \times \cdots \times [\delta, \delta + t_p])^m \\
&\quad = \int_{(\mathbb{R}^d)^m} dx_1 \cdots dx_m \prod_{j=1}^p \int_{0 \le s_1 < \cdots < s_m \le t_j} ds_1 \cdots ds_m \\
&\qquad \times \sum_{\sigma \in \Sigma_m} p_{\delta + s_1}(x_{\sigma(1)}) \prod_{k=2}^m p_{s_k - s_{k-1}}(x_{\sigma(k)} - x_{\sigma(k-1)}).
\end{aligned}
$$

In view of (4.3), therefore,

$$
\begin{aligned}
&\mathbb{E}\alpha([\delta, \delta + \tau_1] \times \cdots \times [\delta, \delta + \tau_p])^m \\
&\quad \le A_\theta (e^{p\delta} \lambda^{-(p+d)})^m (pm)! \theta^{-pm} \qquad \forall\, \theta < p \cdot M^{-(2p - d(p-1))/(2p)}.
\end{aligned} \tag{4.4}
$$

For each $0 \le j \le p$ and the integers $1 \le k_1 < \cdots < k_j \le p$, let $A_{k_1 \cdots k_j}$ be the product of the $p$ sets, of which the $k_1, \dots, k_j$ factors are $[\delta, \delta + \tau_{k_1}], \dots, [\delta, \delta + \tau_{k_j}]$, respectively, and the rest are $[0, \delta]$. Then

$$
\alpha([0, \delta + \tau_1] \times \cdots \times [0, \delta + \tau_p]) = \sum_{j=0}^p \sum_{1 \le k_1 < \cdots < k_j \le p} \alpha(A_{k_1 \cdots k_j}).
$$

By the fact that

$$
\alpha([0, \tau_1] \times \cdots \times [0, \tau_p]) \le \alpha([0, \delta + \tau_1] \times \cdots \times [0, \delta + \tau_p])
$$

and by Hölder's triange inequality,

$$
[\mathbb{E}\alpha([0, \tau_1] \times \cdots \times [0, \tau_p])^m]^{1/m}
$$

$$
\begin{aligned}
&\leq \sum_{j=0}^{p} \sum_{1\leq k_1<\cdots<k_j\leq p} [\mathbb{E}\alpha(A_{k_1\cdots k_j})^m]^{1/m} \\
&\leq [\mathbb{E}\alpha([\delta,\delta+\tau_1]\times\cdots\times[\delta,\delta+\tau_p])^m]^{1/m} \\
&\quad + \sum_{j=0}^{p-1} \binom{p}{j} [\mathbb{E}\alpha([\delta,\delta+\tau]^p)^m]^{j/mp} [\mathbb{E}\alpha([0,\delta]^p)^m]^{(p-j)/(mp)},
\end{aligned}
\tag{4.5}
$$

where the second step follows from Lemma 4.1.

Let $\tau_1',\dots,\tau_p'$ be independent exponential times with parameter $p^{-1}$. Similar to (4.4), there is a constant $C>0$ such that

$$
\mathbb{E}\alpha([\delta,\delta+\tau_1']\times\cdots\times[\delta,\delta+\tau_p'])^m \leq (pm)!C^m, \qquad m\geq 1.
$$

Notice that $\tau \overset{d}{=} \min\{\tau_1',\dots,\tau_p'\}$,

$$
\mathbb{E}\alpha([\delta,\delta+\tau]^p)^m \leq \alpha([\delta,\delta+\tau_1']\times\cdots\times[\delta,\delta+\tau_p'])^m \leq (pm)!C^m, \qquad m\geq 1.
$$

Taking (1.12), (4.4) and (4.5) into account, we obtain

$$
\begin{aligned}
&\mathbb{E}\alpha([0,\tau_1]\times\cdots\times[0,\tau_p])^m \\
&\quad \leq A_\theta(1+o(1))^m (e^{p\delta}\lambda^{-(p+d)})^m (pm)!\theta^{-pm} \qquad (m\to\infty).
\end{aligned}
\tag{4.6}
$$

From (1.10) we have

$$
\begin{aligned}
\mathbb{E}\alpha([0,1]^p)^m \leq{}& A_\theta(1+o(1))^m (e^{p\delta}\lambda^{-(p+d)})^m p^{1+((2p-d(p-1))/2)m} \\
&\times \Gamma\biggl(1+\frac{2p-d(p-1)}{2}m\biggr)^{-1}(pm)!\theta^{-pm}.
\end{aligned}
$$

An estimate via the Stirling formula gives

$$
\begin{aligned}
&\limsup_{k\to\infty} \sqrt[k]{\frac{\mathbb{E}\alpha([0,1]^p)^{2k/(d(p-1))}}{k!}} \\
&\quad \leq (e^{p\delta}\lambda^{-(p+d)})^{2/(d(p-1))} \frac{2}{d(p-1)} \\
&\qquad \times \biggl(\frac{p}{\theta}\biggr)^{2p/(d(p-1))} \biggl(\frac{2p}{2p-d(p-1)}\biggr)^{(2p-d(p-1))/(d(p-1))}.
\end{aligned}
$$

Take the limit in the order $\lambda\to 1$, $\delta\to 0$ and $\theta\to pM^{-(2p-d(p-1))/(2p)}$. So we have

$$
\limsup_{k\to\infty} \sqrt[k]{\frac{\mathbb{E}\alpha([0,1]^p)^{2k/(d(p-1))}}{k!}}
$$

$$\le \frac{2}{d(p-1)}\left(\frac{2pM}{2p-d(p-1)}\right)^{(2p-d(p-1))/(d(p-1))}$$
$$= \frac{2}{p}\kappa(d,p)^{4p/(d(p-1))},$$

where the equality follows from Lemma A.2.

Therefore, by Taylor's expansion we can see that

$$\mathbb{E}\exp\{\gamma\alpha([0,1]^p)^{2/(d(p-1))}\} < \infty \qquad \forall\,\gamma < \frac{p}{2}\kappa(d,p)^{-4p/(d(p-1))}. \tag{4.7}$$

Finally, (4.1) follows from Chebyshev's inequality.

REMARK 4.1. To our surprise, the estimate given in (1.10) turns out to be sharp enough to maintain the right constants. We point out a possible connection between our results and those established by König and Mörters (2002). From (4.6), we can see that

$$\limsup_{m\to\infty}\frac{1}{m}\log\frac{1}{(m!)^p}\mathbb{E}\alpha([0,\tau_1]\times\cdots\times[0,\tau_p])^m$$
$$\le \frac{2p-d(p-1)}{2}\log(2p-d(p-1))$$
$$+\frac{d(p-1)}{2}\log d(p-1) + p\log\frac{\kappa(d,p)^2}{p}.$$

From the relationship (1.10) and our main result, Theorem 2.1, the opposite relationship follows. So we have

$$\begin{aligned}\lim_{m\to\infty}\frac{1}{m}\log\frac{1}{(m!)^p}\mathbb{E}\alpha([0,\tau_1]\times\cdots\times[0,\tau_p])^m &\\ = \frac{2p-d(p-1)}{2}\log(2p-d(p-1)) &\\ +\frac{d(p-1)}{2}\log d(p-1) + p\log\frac{\kappa(d,p)^2}{p}.&\end{aligned} \tag{4.8}$$

This result takes a form close to that given in Proposition 2.2 of König and Mörters (2002), which provides the key estimate for their theorems. In the paper by König and Mörters, the exponential times are replaced by exit times and the intersection local time is limited to a bounded domain. It may be of some interest in the future study to understand how exactly they are related to each other.

We close this section with the following lemma.

LEMMA 4.1. *For any bounded Borel sets $A_1, \dots, A_p \in \mathbb{R}^+$ and for any integer $m \ge 1$,*

$$\mathbb{E}[\alpha(A_1 \times \cdots \times A_p)^m] \le \prod_{j=1}^{p} (\mathbb{E}[\alpha(A_j^p)^m])^{1/p}.$$

PROOF. Write

$$(A_j)_<^m = \{(s_1, \dots, s_m) \in (A_j)^m;\ s_1 < \cdots < s_m\}.$$

Then our lemma follows from Le Gall's moment formula [Le Gall (1990), Theorem 1, page 182] and Hölder's inequality [recall our convention $\sigma(0) \equiv 0$ and $s_0 \equiv 0$]:

$$\begin{aligned}
&\mathbb{E}[\alpha(A_1 \times \cdots \times A_p)^m] \\
&\quad = \int_{(\mathbb{R}^d)^m} dx_1 \cdots dx_m \prod_{j=1}^{p} \biggl[\int_{(A_j)_<^m} ds_1 \cdots ds_m \\
&\qquad\qquad \times \sum_{\sigma \in \Sigma_m} \prod_{k=1}^{m} p_{s_k - s_{k-1}}(x_{\sigma(k)} - x_{\sigma(k-1)})\biggr] \\
&\quad \le \prod_{j=1}^{p} \biggl\{\int_{(\mathbb{R}^d)^m} dx_1 \cdots dx_m \biggl[\int_{(A_j)_<^m} ds_1 \cdots ds_m \\
&\qquad\qquad \times \sum_{\sigma \in \Sigma_m} \prod_{k=1}^{m} p_{s_k - s_{k-1}}(x_{\sigma(k)} - x_{\sigma(k-1)})\biggr]^p\biggr\}^{1/p} \\
&\quad = \prod_{j=1}^{p} (\mathbb{E}[\alpha(A_j^p)^m])^{1/p}. \qquad \square
\end{aligned}$$

**5. Upper bound in Theorem 2.2.** The approach used in Section 4 is no longer applicable to the case of random walks, mainly due to the absence of a scaling property. To begin with, we introduce the following inequality which is of interest for its own sake.

THEOREM 5.1. *Given positive integers $n_1, \dots, n_a$ and $m$,*

$$(\mathbb{E}I_{n_1+\cdots+n_a}^m)^{1/p} \le \sum_{\substack{k_1+\cdots+k_a=m \\ k_1,\dots,k_a \ge 0}} \frac{m!}{k_1! \cdots k_a!} (\mathbb{E}I_{n_1}^{k_1})^{1/p} \cdots (\mathbb{E}I_{n_a}^{k_a})^{1/p}. \tag{5.1}$$

*Consequently, for any $\lambda > 0$,*

$$\sum_{m=0}^{\infty} \frac{\lambda^m}{m!} (\mathbb{E}I_{n_1+\cdots+n_a}^m)^{1/p} \le \prod_{i=1}^{a} \sum_{m=0}^{\infty} \frac{\lambda^m}{m!} (\mathbb{E}I_{n_i}^m)^{1/p}. \tag{5.2}$$

PROOF. Let $l(n,x)$ be the local time generated by $S(n)$,

$$l(n,x) = \sum_{k=1}^{n} \mathbb{1}_{\{S(k)=x\}}, \qquad n = 1, 2, \ldots,$$

and let $l_1(n,x), \ldots, l_p(n,x)$ be the local times of the independent random walks $\{S_1(n)\}, \ldots, \{S_p(n)\}$, respectively. Then

$$I_n = \sum_{x\in\mathbb{Z}^d} \prod_{j=1}^{p} l_j(n,x). \tag{5.3}$$

Write $n_0 = 0$ and

$$\Delta_i(x) = l(n_0 + \cdots + n_i, x) - l(n_0 + \cdots + n_{i-1}, x), \qquad x \in \mathbb{Z}^d; i = 1, \ldots, a,$$

$$\Delta_{ij}(x) = l_j(n_0 + \cdots + n_i, x) - l_j(n_0 + \cdots + n_{i-1}, x), \qquad x \in \mathbb{Z}^d; i = 1, \ldots, a; j = 1, \ldots, p.$$

Then

$$\begin{aligned}(\mathbb{E}I^m_{n_1+\cdots+n_a})^{1/p} &= \left(\sum_{x_1,\ldots,x_m} \left[\mathbb{E}\prod_{k=1}^{m}\sum_{i=1}^{a}\Delta_i(x_k)\right]^p\right)^{1/p} \\ &= \left(\sum_{x_1,\ldots,x_m}\left[\sum_{i_1,\ldots,i_m=1}^{a} \mathbb{E}(\Delta_{i_1}(x_1)\cdots\Delta_{i_m}(x_m))\right]^p\right)^{1/p} \\ &\le \sum_{i_1,\ldots,i_m=1}^{a}\left(\sum_{x_1,\ldots,x_m}[\mathbb{E}(\Delta_{i_1}(x_1)\cdots\Delta_{i_m}(x_m))]^p\right)^{1/p}.\end{aligned}$$

Given integers $i_1, \ldots, i_m$ between 1 and $a$, let $k_1, \ldots, k_a$ be the numbers of $1, \ldots, a$, respectively, among $i_1, \ldots, i_m$. Then $k_1 + \cdots + k_a = m$. To prove our conclusion, it suffices to show that

$$\sum_{x_1,\ldots,x_m}[\mathbb{E}(\Delta_{i_1}(x_1)\cdots\Delta_{i_m}(x_m))]^p \le \mathbb{E}I^{k_1}_{n_1}\cdots\mathbb{E}I^{k_a}_{n_a}.$$

Without loss of generality, we may only consider the case when $k_1, \ldots, k_a \ge 1$. Under the notation $\bar{x}_i = (x^i_1, \ldots, x^i_{k_i}) \in (\mathbb{Z}^d)^{k_i}$, we set

$$\phi_i(\bar{x}_i) = \mathbb{E}(l(n_i, x^i_1)\cdots l(n_i, x^i_{k_i})).$$

It is easy to see that

$$\sum_{\bar{x}_i}\phi^p_i(\bar{x}_i) = \mathbb{E}I^{k_i}_{n_i}, \qquad i = 1, \ldots, a.$$

Define

$$\bar{S}^i(k) = (\overbrace{S(k), \dots, S(k)}^{k_i}) \quad \text{and} \quad \bar{S}^i_j(k) = (\overbrace{S_j(k), \dots, S_j(k)}^{k_i}), \qquad k = 1, 2, \dots,$$

where $1 \le i \le a$ and $1 \le j \le p$. Then

$$\begin{aligned}
&\sum_{x_1,\dots,x_m} [\mathbb{E}(\Delta_{i_1}(x_1)\cdots\Delta_{i_m}(x_m))]^p \\
&\quad = \sum_{\bar{x}_1}\cdots\sum_{\bar{x}_a}\left[\mathbb{E}\prod_{i=1}^a \Delta_i(x_1^i)\cdots\Delta_i(x_{k_i}^i)\right]^p \\
&\quad = \sum_{\bar{x}_1}\cdots\sum_{\bar{x}_a}\left[\mathbb{E}\left\{\left(\prod_{i=1}^{a-1}\Delta_i(x_1^i)\cdots\Delta_i(x_{k_i}^i)\right)\phi_a(\bar{x}_a - \bar{S}^a(n-n_a))\right\}\right]^p.
\end{aligned}$$

Notice that

$$\begin{aligned}
&\sum_{\bar{x}_a}\left[\mathbb{E}\left\{\left(\prod_{i=1}^{a-1}\Delta_i(x_1^i)\cdots\Delta_i(x_{k_i}^i)\right)\phi_a(\bar{x}_a - \bar{S}^a(n-n_a))\right\}\right]^p \\
&\quad = \sum_{\bar{x}_a}\mathbb{E}\left\{\prod_{j=1}^p\left(\prod_{i=1}^{a-1}\Delta_{ij}(x_1^i)\cdots\Delta_{ij}(x_{k_i}^i)\right)\phi_a(\bar{x}_a - \bar{S}^a_j(n-n_a))\right\} \\
&\quad = \mathbb{E}\left\{\left(\prod_{j=1}^p\prod_{i=1}^{a-1}\Delta_{ij}(x_1^i)\cdots\Delta_{ij}(x_{k_i}^i)\right)\sum_{\bar{x}_a}\prod_{j=1}^p\phi_a(\bar{x}_a - \bar{S}^a_j(n-n_a))\right\} \\
&\quad \le \mathbb{E}\left\{\left(\prod_{j=1}^p\prod_{i=1}^{a-1}\Delta_{ij}(x_1^i)\cdots\Delta_{ij}(x_{k_i}^i)\right)\prod_{j=1}^p\left(\sum_{\bar{x}_a}\phi_a^p(\bar{x}_a - \bar{S}^a_j(n-n_a))\right)^{1/p}\right\} \\
&\quad = \mathbb{E}\left\{\left(\prod_{j=1}^p\prod_{i=1}^{a-1}\Delta_{ij}(x_1^i)\cdots\Delta_{ij}(x_{k_i}^i)\right)\sum_{\bar{x}_a}\phi_a^p(\bar{x}_a)\right\} \\
&\quad = \left(\mathbb{E}\prod_{i=1}^{a-1}\Delta_i(x_1^i)\cdots\Delta_i(x_{k_i}^i)\right)^p\cdot\mathbb{E}I_{n_a}^{k_a}.
\end{aligned}$$

So we have

$$\begin{aligned}
&\sum_{\bar{x}_1}\cdots\sum_{\bar{x}_a}\left[\mathbb{E}\prod_{i=1}^a\Delta_i(x_1^i)\cdots\Delta_i(x_{k_i}^i)\right]^p \\
&\qquad \le \mathbb{E}I_{n_a}^{k_a}\cdot\sum_{\bar{x}_1}\cdots\sum_{\bar{x}_{a-1}}\left[\mathbb{E}\prod_{i=1}^{a-1}\Delta_i(x_1^i)\cdots\Delta_i(x_{k_i}^i)\right]^p.
\end{aligned}$$

Repeat this procedure:

$$\sum_{x_1,\dots,x_m}[\mathbb{E}(\Delta_{i_1}(x_1)\cdots\Delta_{i_m}(x_m))]^p \le \mathbb{E}I_{n_1}^{k_1}\cdots\mathbb{E}I_{n_a}^{k_a}.$$

The second half of Theorem 5.1 follows from the computation

$$\begin{aligned}
&\sum_{m=0}^{\infty}\frac{\lambda^m}{m!}(\mathbb{E}I^m_{n_1+\cdots+n_a})^{1/p}\\
&\quad\le\sum_{m=0}^{\infty}\lambda^m\sum_{\substack{k_1+\cdots+k_a=m\\k_1,\ldots,k_a\ge0}}\frac{1}{k_1!\cdots k_a!}(\mathbb{E}I^{k_1}_{n_1})^{1/p}\cdots(\mathbb{E}I^{k_a}_{n_a})^{1/p}\\
&\quad=\sum_{k_1,\ldots,k_a=1}^{\infty}\lambda^{k_1+\cdots+k_a}\frac{1}{k_1!\cdots k_a!}(\mathbb{E}I^{k_1}_{n_1})^{1/p}\cdots(\mathbb{E}I^{k_a}_{n_a})^{1/p}\\
&\quad=\prod_{i=1}^{a}\sum_{m=0}^{\infty}\frac{\lambda^m}{m!}(\mathbb{E}I^m_{n_i})^{1/p}.
\end{aligned}$$

□

Theorem 5.1 applies to our situation in two different ways. The first application is given in the following lemma.

LEMMA 5.2. *For given $\lambda>0$, there is a positive sequence $\{C_m\}_{m\ge0}$ such that*

$$\sup_n(n^{-((2p-d(p-1))/2)m}\mathbb{E}I^m_n)^{1/p}\le C_m,\qquad m\ge0,$$

*and*

$$\sum_{m=0}^{\infty}\frac{\lambda^m}{m!}C_m<\infty.$$

PROOF. Write $P^k(x)=\mathbb{P}\{S_k=x\}$ $(k\ge1,\ x\in\mathbb{Z}^d)$ and $P^0(x)=\delta_0(x)$. For any $m\ge1$, let $\Sigma_m$ be the set of the permutations of $\{1,\ldots,m\}$. Under the convention $\sigma(0)=0$ and $i_0=0$,

$$\begin{aligned}
\mathbb{E}I^m_n&=\sum_{x_1,\ldots,x_m}\left[\sum_{\sigma\in\Sigma_m}\sum_{1\le i_1\le\cdots\le i_m\le n}\prod_{k=1}^{m}P^{i_k-i_{k-1}}(x_{\sigma(k)}-x_{\sigma(k-1)})\right]^p\\
&\le(m!)^p\sum_{x_1,\ldots,x_m}\left[\sum_{1\le i_1\le\cdots\le i_m\le n}\prod_{k=1}^{m}P^{i_k-i_{k-1}}(x_k)\right]^p\\
&\le(m!)^p\sum_{x_1,\ldots,x_m}\left[\prod_{k=1}^{m}\sum_{i=0}^{n}P^i(x_k)\right]^p\\
&=(m!)^p\left[\sum_x\left(\sum_{i=0}^{n}P^i(x)\right)^p\right]^m.
\end{aligned}$$

Notice that

$$\sum_x \left(\sum_{i=0}^n P^i(x)\right)^p \sim \sum_x \left(\sum_{i=1}^n P^i(x)\right)^p = \mathbb{E}I_n, \qquad (n\to\infty).$$

In connection with the weak law given in (1.4), we have

$$\mathbb{E}I_n = O(n^{(2p-d(p-1))/2}), \qquad (n\to\infty).$$

Therefore, there is a $C>0$ such that

$$\sup_n (n^{-((2p-d(p-1))/2)m}\mathbb{E}I_n^m)^{1/p} \le m!C^m, \qquad m=0,1,\ldots.$$

When $\lambda C<1$, the lemma follows if we take $C_m = m!C^m$. In the case $\lambda C\ge 1$, we choose a small $\delta>0$ such that

$$\lambda\delta^{(2p-d(p-1))/(2p)}C<1.$$

Let $a=[\delta^{-1}]+1$. Then $n\le a[\delta n]$. By (5.1) in Theorem 5.1,

$$\begin{aligned}
(\mathbb{E}I_n^m)^{1/p} &\le \sum_{k_1+\cdots+k_a=m} \frac{m!}{k_1!\cdots k_a!}(\mathbb{E}I_{[\delta n]}^{k_1})^{1/p}\cdots(\mathbb{E}I_{[\delta n]}^{k_a})^{1/p} \\
&\le \sum_{k_1+\cdots+k_a=m} \frac{m!}{k_1!\cdots k_a!}\prod_{i=1}^a k_i!C^{k_i}[\delta n]^{((2p-d(p-1))/(2p))k_i} \\
&= m!C^m[\delta n]^{((2p-d(p-1))/(2p))m}\binom{m+a-1}{m} \\
&\le m!(\delta^{(2p-d(p-1))/(2p)}C)^m\binom{m+a-1}{m}n^{((2p-d(p-1))/(2p))m} \qquad \forall\, n,
\end{aligned}$$

where the equality follows from the fact that the equation $k_1+\cdots+k_a=m$ has $\binom{m+a-1}{m}$ nonnegative integer solutions.

Therefore, the desired conclusion follows with

$$C_m = m!(\delta^{((2p-d(p-1))/(2p))}C)^m\binom{m+a-1}{m}, \qquad m=0,1,\ldots. \qquad \square$$

LEMMA 5.3. *Let $\{Z_\varepsilon\}$ be a family of nonnegative random variables and let $p\ge 1$ be an integer. Let $I(\lambda)$ be a lower semicontinuous, nondecreasing, nonnegative function on $[0,\infty)$ such that $I(0)=0$, $I(|\cdot|^p)$ is convex on $(-\infty,\infty)$ and that $I(\lambda)\to\infty$ as $\lambda\to\infty$.*

(i) *If*

$$\limsup_{\varepsilon\to 0^+}\varepsilon\log\mathbb{P}\{Z_\varepsilon\ge\lambda\}\le -I(\lambda), \qquad (\lambda>0) \tag{5.4}$$

*and if* $b>0$ *satisfies*

$$\lim_{N\to\infty}\limsup_{\varepsilon\to 0^+}\varepsilon\log\sum_{m=1}^{\infty}\frac{(b\varepsilon^{-1})^m}{m!}(\mathbb{E}Z_\varepsilon^m \mathbb{1}_{\{Z_\varepsilon\ge N\}})^{1/p}=-\infty, \tag{5.5}$$

*then*

$$\limsup_{\varepsilon\to 0^+}\varepsilon\log\left(\sum_{m=0}^{\infty}\frac{(b\varepsilon^{-1})^m}{m!}(\mathbb{E}Z_\varepsilon^m)^{1/p}\right)\le\sup_{\lambda>0}\{b\lambda^{1/p}-p^{-1}I(\lambda)\}. \tag{5.6}$$

(ii) *Conversely, if* (5.6) *is satisfied for all* $b>0$, *then* (5.4) *holds for all* $\lambda>0$.

(iii) *In addition, the condition* (5.5) *is satisfied if there is a* $b'>2pb$ *such that*

$$\limsup_{\varepsilon\to 0^+}\varepsilon\log\mathbb{E}\exp\{\varepsilon^{-1}b'Z_\varepsilon^{1/p}\}<\infty. \tag{5.7}$$

Remark 5.1. By the convention used in the area of large deviations, the notation $\varepsilon$ may also be used for a positive sequence approaching zero, in which case $Z_\varepsilon$ represents a random sequence.

Proof of Lemma 5.1. Part (i) follows from an argument almost identical to that for Lemma 4.3.6 in Dembo and Zeitouni (1998).

To prove (ii), write

$$\Psi(b)=\sup_{\lambda\in\mathbb{R}}\{b\lambda-p^{-1}I(|\lambda|^p)\},\qquad b\in\mathbb{R}.$$

If $b>0$,

$$\Psi(b)=\sup_{\lambda>0}\{b\lambda^{1/p}-p^{-1}I(\lambda)\}.$$

By Chebyshev's inequality, as $b>0$,

$$\lambda^{m/p}(b\varepsilon^{-1})^m(\mathbb{P}\{Z_\varepsilon\ge\lambda\})^{1/p}\le(b\varepsilon^{-1})^m(\mathbb{E}Z_\varepsilon^m)^{1/p}$$

for any integer $m\ge 0$. Summing up gives

$$\exp(b\lambda^{1/p}\varepsilon^{-1})(\mathbb{P}\{Z_\varepsilon\ge\lambda\})^{1/p}\le\sum_{m=0}^{\infty}\frac{(b\varepsilon^{-1})^m}{m!}(\mathbb{E}Z_\varepsilon^m)^{1/p}.$$

Hence

$$\limsup_{\varepsilon\to 0^+}\varepsilon\log\mathbb{P}\{Z_\varepsilon\ge\lambda\}\le -p\{\lambda^{1/p}b-\Psi(b)\}. \tag{5.8}$$

Notice that

$$\sup_{b>0}\{\lambda^{1/p}b-\Psi(b)\}=\sup_{b\in\mathbb{R}}\{\lambda^{1/p}b-\Psi(b)\}=p^{-1}I(\lambda)\qquad(\lambda>0),$$

where the first equality follows from the fact

$$\sup_{b>0}\{\lambda^{1/p}b-\Psi(b)\}\ge 0 \quad\text{and}\quad \sup_{b\le 0}\{\lambda^{1/p}b-\Psi(b)\}\le 0$$

and the second follows from the duality lemma [see, e.g., Dembo and Zeitouni (1998), Lemma 4.5.8]. Taking the supremum over $b>0$ on the right-hand side of (5.8) gives (5.4).

We now prove (iii). From the relationship $(2pm)!\le(2p)^{2pm}(m!)^{2p}$,

$$\begin{aligned}
&\sum_{m=0}^{\infty}\frac{(b\varepsilon^{-1})^m}{m!}(\mathbb{E}Z_\varepsilon^m\mathbb{1}_{\{Z_\varepsilon\ge N\}})^{1/p}\\
&\qquad\le(\mathbb{P}\{Z_\varepsilon\ge N\})^{1/(2p)}\sum_{m=0}^{\infty}\frac{(b\varepsilon^{-1})^m}{m!}(\mathbb{E}Z_\varepsilon^{2m})^{1/(2p)}\\
&\qquad\le(\mathbb{P}\{Z_\varepsilon\ge N\})^{1/(2p)}\sum_{m=0}^{\infty}\left(\frac{2pb}{b'}\right)^m\left(\frac{(\varepsilon^{-1}b')^{2mp}\mathbb{E}Z_\varepsilon^{2m}}{(2mp)!}\right)^{1/(2p)}\\
&\qquad\le\left(1-\frac{2pb}{b'}\right)^{-1}(\mathbb{P}\{Z_\varepsilon\ge N\})^{1/(2p)}(\mathbb{E}\exp\{\varepsilon^{-1}b'Z_\varepsilon^{1/p}\})^{1/(2p)}.
\end{aligned}$$

Hence

$$\begin{aligned}
&\limsup_{\varepsilon\to0^+}\varepsilon\log\sum_{m=1}^{\infty}\frac{(b\varepsilon^{-1})^m}{m!}(\mathbb{E}Z_\varepsilon^m\mathbb{1}_{\{Z_\varepsilon\ge N\}})^{1/p}\\
&\qquad\le\frac{1}{2p}\left\{-I(N)+\limsup_{\varepsilon\to0^+}\varepsilon\log\mathbb{E}\exp\{\varepsilon^{-1}b'Z_\varepsilon^{1/p}\}\right\}.
\end{aligned}$$

Letting $N\to\infty$ gives (5.5). $\square$

We are ready to prove the upper bound for Theorem 2.2:

$$\begin{aligned}
&\limsup_{n\to\infty}\frac{1}{b_n}\log\mathbb{P}\{I_n\ge\lambda n^{(2p-d(p-1))/2}b_n^{(d(p-1))/2}\}\\
&\qquad\le-\frac{p}{2}\det(\Gamma)^{1/d}\kappa(d,p)^{-4p/(d(p-1))}\lambda^{2/(d(p-1))}.
\end{aligned}\tag{5.9}$$

The proof is the second application of Theorem 5.1. Let $t>0$ be fixed, and let $t_n=[tn/b_n]$ and $\gamma_n=[n/t_n]$. Then $n\le t_n(\gamma_n+1)$. By (5.2) in Theorem 5.1,

$$\begin{aligned}
&\sum_{m=0}^{\infty}\frac{1}{m!}\theta^m\left(\frac{b_n}{n}\right)^{((2p-d(p-1))/(2p))m}(\mathbb{E}I_n^m)^{1/p}\\
&\qquad\le\left(\sum_{m=0}^{\infty}\frac{1}{m!}\theta^m\left(\frac{b_n}{n}\right)^{((2p-d(p-1))/(2p))m}(\mathbb{E}I_{t_n}^m)^{1/p}\right)^{\gamma_n+1}
\end{aligned}$$

for any $\theta > 0$. In view of Lemma 5.2, by the weak convergence given in (1.4) and the dominated convergence theorem,

$$\sum_{m=0}^{\infty} \frac{1}{m!} \theta^m \left(\frac{b_n}{n}\right)^{((2p-d(p-1))/(2p))m} (\mathbb{E} I_{t_n}^m)^{1/p}$$
$$\to \sum_{m=0}^{\infty} \frac{1}{m!} \theta^m t^{((2p-d(p-1))/(2p))m} \det(\Gamma)^{-((p-1)/(2p))m} (\mathbb{E}\alpha([0,1]^p)^m)^{1/p}$$

as $n \to \infty$. Hence,

$$\limsup_{n\to\infty} \frac{1}{b_n} \log\left(\sum_{m=0}^{\infty} \frac{1}{m!} \theta^m \left(\frac{b_n}{n}\right)^{((2p-d(p-1))/(2p))m} (\mathbb{E} I_n^m)^{1/p}\right)$$
$$\le \frac{1}{t} \log\left(\sum_{m=0}^{\infty} \frac{1}{m!} \theta^m t^{((2p-d(p-1))/(2p))m} \times \det(\Gamma)^{-((p-1)/(2p))m} (\mathbb{E}\alpha([0,1]^p)^m)^{1/p}\right). \tag{5.10}$$

In addition, (4.1) implies that

$$\limsup_{t\to\infty} \frac{1}{t} \log \mathbb{E} \exp\{b' t^{(2p-d(p-1))/(2p)} \alpha([0,1]^p)^{1/p}\} < \infty$$

for any $b' > 0$. Hence condition (5.7) is satisfied with $\varepsilon = t^{-1}$ and with

$$Z_\varepsilon = t^{-(d(p-1))/2} \alpha([0,1]^p),$$

and [by (4.1)] condition (5.4) is satisfied with

$$I(\lambda) = \frac{p}{2} \kappa(d,p)^{-4p/(d(p-1))} \lambda^{2/(d(p-1))}.$$

According to (i) of Lemma 5.3,

$$\limsup_{t\to\infty} \frac{1}{t} \log\left(\sum_{m=0}^{\infty} \frac{1}{m!} \theta^m t^{((2p-d(p-1))/(2p))m} \times \det(\Gamma)^{-((p-1)/(2p))m} (\mathbb{E}\alpha([0,1]^p)^m)^{1/p}\right)$$
$$\le \sup_{\lambda>0} \left\{\det(\Gamma)^{-(p-1)/(2p)} \theta \lambda^{1/p} - \frac{1}{p} \frac{p}{2} \kappa(d,p)^{-4p/(d(p-1))} \lambda^{2/(d(p-1))}\right\}$$
$$= \sup_{\lambda>0} \left\{\theta \lambda^{1/p} - \frac{1}{p} \frac{p}{2} \det(\Gamma)^{1/d} \kappa(d,p)^{-4p/(d(p-1))} \lambda^{2/(d(p-1))}\right\}.$$

Letting $t \to \infty$ in (5.10) gives

$$\begin{aligned}(5.11)\quad &\limsup_{n\to\infty} \frac{1}{b_n} \log\left(\sum_{m=0}^{\infty} \frac{1}{m!}\theta^m \left(\frac{b_n}{n}\right)^{((2p-d(p-1))/(2p))m} (\mathbb{E}I_n^m)^{1/p}\right) \\ &\qquad \le \sup_{\lambda>0}\left\{\theta\lambda^{1/p} - \frac{1}{p}\frac{p}{2}\det(\Gamma)^{1/d}\kappa(d,p)^{-4p/(d(p-1))}\lambda^{2/(d(p-1))}\right\}.\end{aligned}$$

Finally, the desired upper bound (5.9) follows from (ii) of Lemma 5.3 with $\varepsilon = b_n^{-1}$ and $Z_\varepsilon$ replaced by

$$n^{-(2p-d(p-1))/2} b_n^{-(d(p-1))/2} I_n.$$

**6. Lower bounds.** For given $\bar{x} = (x_1^o, \dots, x_p^o) \in (\mathbb{R}^d)^p$, we introduce the notation $\mathbb{P}^{\bar{x}}$ for the probability induced by independent $d$-dimensional Brownian motions $W_1, \dots, W_p$ starting at $x_1^o, \dots, x_p^o$, respectively. Without causing confusion, for given $\bar{x} = (x_1^o, \dots, x_p^o) \in (\mathbb{Z}^d)^p$, we also use $\mathbb{P}^{\bar{x}}$ for the probability induced by the random walks $S_1(n), \dots, S_p(n)$ in the case when $S_1(n), \dots, S_p(n)$ start at $x_1^o, \dots, x_p^o$, respectively. The notation $\mathbb{E}^{\bar{x}}$ denotes the expectation that corresponds to $\mathbb{P}^{\bar{x}}$. To be consistent with the notation we used before, we have $\mathbb{P}^{(0,\dots,0)} = \mathbb{P}$ and $\mathbb{E}^{(0,\dots,0)} = \mathbb{E}$. Write

$$\|\bar{x}\| = \max_{1\le j\le p} |x_j^o|, \qquad \bar{x} = (x_1^o, \dots, x_p^o) \in (\mathbb{R}^d)^p.$$

THEOREM 6.1. *For any constant $C > 0$,*

$$(6.1)\quad \liminf_{t\to\infty} \frac{1}{t} \log \inf_{\|\bar{x}\|\le C} \mathbb{P}^{\bar{x}}\{\alpha([0,1]^p) \ge t^{(d(p-1))/2}\} \ge -\frac{p}{2}\kappa(d,p)^{-4p/(d(p-1))}.$$

*Given $\lambda > 0$,*

$$\begin{aligned}(6.2)\quad &\liminf_{n\to\infty} \frac{1}{b_n} \log \inf_{\|\bar{x}\|\le C\sqrt{n}} \mathbb{P}^{\bar{x}}\{I_n \ge \lambda n^{(2p-d(p-1))/2} b_n^{(d(p-1))/2}\} \\ &\qquad \ge -\frac{p}{2}\det(\Gamma)^{1/d}\kappa(d,p)^{-4p/(d(p-1))}\lambda^{2/(d(p-1))}.\end{aligned}$$

PROOF. Due to similarity, we only prove (6.2). We first proceed under the additional assumption that the random walk $\{S(n)\}$ is aperiodic: The greatest common factor of the set $\{n \ge 1; \mathbb{P}\{S(n) = 0\} > 0\}$ is 1. Let $M > 0$ be given as in (4.2). To prove (6.2), it is sufficient to show that for any $\theta > 0$,

$$\begin{aligned}(6.3)\quad &\liminf_{n\to\infty} \frac{1}{b_n} \log \inf_{\|\bar{x}\|\le C\sqrt{n}} \mathbb{E}^{\bar{x}} \exp\left\{\theta\left(\frac{b_n}{n}\right)^{(2p-d(p-1))/(2p)} I_n^{1/p}\right\} \\ &\qquad \ge \theta^{2p/(2p-d(p-1))} p^{-(d(p-1))/(2p-d(p-1))} \det(\Gamma)^{-(p-1)/(2p-d(p-1))} M.\end{aligned}$$

Indeed, since $\mathbb{E}I_n^{m/p} \le (\mathbb{E}I_n^m)^{1/p}$ by (5.11), we have

$$\limsup_{n\to\infty} \frac{1}{b_n} \log \mathbb{E}\exp\biggl\{\theta\biggl(\frac{b_n}{n}\biggr)^{(2p-d(p-1))/(2p)} I_n^{1/p}\biggr\} < \infty \qquad \forall\,\theta > 0.$$

Furthermore, by (5.9) and Lemma 4.3.6 in Dembo and Zeitouni (1998), the above lim sup is bounded by

$$\begin{aligned}&\sup_{\lambda>0}\biggl\{\theta\lambda^{1/p} - \frac{p}{2}\det(\Gamma)^{1/d}\kappa(d,p)^{-4p/(d(p-1))}\lambda^{2/(d(p-1))}\biggr\}\\ &\qquad = \theta^{2p/(2p-d(p-1))}p^{-(d(p-1))/(2p-d(p-1))}\det(\Gamma)^{-(p-1)/(2p-d(p-1))}M,\end{aligned}$$

where the equality partially follows from Lemma A.2. Combining this with (6.3), by the Gärtner–Ellis theorem [Dembo and Zeitouni (1998), Theorem 2.3.6],

$$\begin{aligned}\liminf_{n\to\infty}\frac{1}{b_n}\log&\inf_{\|\bar{x}\|\le C\sqrt{n}}\mathbb{P}^{\bar{x}}\{I_n \ge \lambda n^{(2p-d(p-1))/2}b_n^{(d(p-1))/2}\}\\ &\ge -\sup_{\theta>0}\{\lambda\theta - \theta^{2p/(2p-d(p-1))}p^{-(d(p-1))/(2p-d(p-1))}\\ &\qquad\qquad\qquad \times\det(\Gamma)^{-(p-1)/(2p-d(p-1))}M\}\\ &= -\frac{p}{2}\det(\Gamma)^{1/d}\kappa(d,p)^{-4p/(d(p-1))}\lambda^{2/(d(p-1))},\end{aligned}$$

where again the equality partially follows from Lemma A.2.

We now prove (6.3). Let $\varepsilon > 0$ and $u > 0$ be small numbers and let $\bar{x} = (x_1^o,\dots,x_p^o) \in (\mathbb{Z}^d)^p$ such that $\|\bar{x}\| \le C\sqrt{n}$. Write

$$B_n(x) = \{y; |y-x| \le \varepsilon\sqrt{n/b_n}\}, \qquad x \in \mathbb{Z}^d,$$

and set $B_n = B_n(0)$. Let $0 < u < 1$ be a small number. For any integer $m \ge 1$,

$$\begin{aligned}&\mathbb{E}^{\bar{x}}\biggl(\sum_x\prod_{j=1}^p[l_j(n+[un],x) - l_j([un],x)]\biggr)^m\\ &\quad = \sum_{x_1,\dots,x_m}\prod_{j=1}^p\mathbb{E}\biggl(\prod_{k=1}^m[l_j(n+[un],x_j^o+x_j) - l_j([un],x_j^o+x_k)]\biggr)\\ &\quad = \sum_{x_1,\dots,x_m}\prod_{j=1}^p\sum_{i_1,\dots,i_m=1}^n\mathbb{E}\biggl(\prod_{k=1}^m\mathbb{1}_{\{S([un]+i_k)=x_j^o+x_k\}}\biggr)\\ &\quad \ge \sum_{x_1,\dots,x_m}\prod_{j=1}^p\sum_{i_1,\dots,i_m=1}^n\mathbb{E}\biggl(\prod_{k=1}^m\sum_{y\in B_n(x_j^o)}\mathbb{1}_{\{S([un])=y\}}\end{aligned}$$

$$\times \mathbb{1}_{\{S([un]+i_k)-S([un])=x_j^o-y+x_k\}}\Bigg).$$

Notice that

$$\prod_{k=1}^m \sum_{y\in B_n(x_j^o)} \mathbb{1}_{\{S([un])=y\}} \cdot \mathbb{1}_{\{S([un]+i_k)-S([un])=x_j^o-y+x_k\}}$$

$$= \sum_{y\in B_n(x_j^o)} \mathbb{1}_{\{S([un])=y\}} \cdot \prod_{k=1}^m \mathbb{1}_{\{S([un]+i_k)-S([un])=x_j^o-y+x_k\}}.$$

Let

$$\gamma_n = \min_{1\le j\le p} \inf_{|x_j^o|\le C\sqrt{n}} \inf_{y\in B_n(x_j^o)} \{\mathbb{P}\{S([un])=y\}\}.$$

Then

$$\mathbb{E}\Bigg(\prod_{k=1}^m \sum_{y\in B_n(x_j^o)} \mathbb{1}_{\{S([un])=y\}} \cdot \mathbb{1}_{\{S([un]+i_k)-S([un])=x_j^o-y+x_k\}}\Bigg)$$

$$= \sum_{y\in B_n(x_j^o)} \mathbb{P}\{S([un])=y\} \cdot \mathbb{E}\Bigg[\prod_{k=1}^m \mathbb{1}_{\{S(i_k)=x_j^o-y+x_k\}}\Bigg]$$

$$\ge \gamma_n \sum_{y\in B_n} \mathbb{E}\Bigg[\prod_{k=1}^m \mathbb{1}_{\{S(i_k)=y+x_k\}}\Bigg].$$

Therefore,

$$\mathbb{E}^{\bar{x}}\Bigg(\sum_x \prod_{j=1}^p [l_j([un]+n,x)-l_j([un],x)]\Bigg)^m$$

$$\ge \gamma_n^p \sum_{x_1,\dots,x_m} \Bigg(\sum_{i_1,\dots,i_m=1}^n \sum_{y\in B_n} \mathbb{E}\Bigg[\prod_{k=1}^m \mathbb{1}_{\{S(i_k)=y+x_k\}}\Bigg]\Bigg)^p$$

$$= \gamma_n^p \sum_{x_1,\dots,x_m} \Bigg(\sum_{y\in B_n} \mathbb{E}\prod_{k=1}^m l(n,y+x_k)\Bigg)^p$$

$$= \gamma_n^p \sum_{x_1,\dots,x_m} \sum_{y_1,\dots,y_p\in B_n} \prod_{j=1}^p \mathbb{E}\Bigg[\prod_{k=1}^m l(n,y_j+x_k)\Bigg]$$

$$= \gamma_n^p \sum_{y_1,\dots,y_p\in B_n} \sum_{x_1,\dots,x_m} \mathbb{E}\Bigg[\prod_{k=1}^m \prod_{j=1}^p l_j(n,y_j+x_k)\Bigg]$$

$$= \gamma_n^p \sum_{y_1,\ldots,y_p \in B_n} \mathbb{E}\biggl(\sum_x \prod_{j=1}^p l_j(n, y_j + x)\biggr)^m.$$

By Jensen's inequality,

$$\begin{aligned}
&\frac{1}{\#\{B_n\}^p} \sum_{y_1,\ldots,y_p \in B_n} \mathbb{E}\biggl(\sum_x \prod_{j=1}^p l_j(n, y_j + x)\biggr)^m \\
&\quad \geq \mathbb{E}\biggl(\frac{1}{\#\{B_n\}^p} \sum_{y_1,\ldots,y_p \in B_n} \sum_x \prod_{j=1}^p l_j(n, y_j + x)\biggr)^m \\
&\quad = \mathbb{E}\biggl(\sum_x \prod_{j=1}^p l_j(n, x, \varepsilon)\biggr)^m,
\end{aligned}$$

where $l_j(n, x, \varepsilon)$ $(1 \leq j \leq p)$ are the same as in Section 3. By (5.3), therefore,

$$\begin{aligned}
\mathbb{E}^{\bar{x}}(I_{n+[un]})^m &\geq \mathbb{E}^{\bar{x}}\biggl(\sum_x \prod_{j=1}^p [l_j([un] + n, x) - l_j([un], x)]\biggr)^m \\
&\geq (\gamma_n \#\{B_n\})^p \mathbb{E}\biggl(\sum_x \prod_{j=1}^p l_j(n, x, \varepsilon)\biggr)^m.
\end{aligned}$$

According to the Remark on page 661 of Le Gall and Rosen (1991), the aperiodicity of the random walk implies

$$\sup_{x \in \mathbb{Z}^d} \biggl| [un]^{d/2} \mathbb{P}\{S([un]) = x\} - \frac{1}{(2\pi)^{d/2} \det(\Gamma)^{1/2}} \exp\biggl\{-\frac{1}{2[un]} \langle x, \Gamma^{-1} x\rangle\biggr\}\biggr| \to 0$$

as $n \to \infty$. Hence

$$\begin{aligned}
\gamma_n &= \frac{1}{[un]^{d/2}} \min_{1 \leq j \leq p} \inf_{|x_j^o| \leq C\sqrt{n}} \inf_{y \in B_n(x_j^o)} \biggl[\exp\biggl\{-\frac{1}{2[un]} \langle y, \Gamma^{-1} y\rangle\biggr\} + o(1)\biggr] \\
&\geq \frac{1}{[un]^{d/2}} \exp\biggl\{-\frac{C^2 + o(1)}{2u\lambda}\biggr\},
\end{aligned}$$

where $\lambda > 0$ is the smallest eigenvalue of $\Gamma$. Notice that $\#\{B_n\} \sim C_d \varepsilon^d (n/b_n)^{d/2}$ as $n \to \infty$. There is a $\delta = \delta(\varepsilon) > 0$, such that for any integer $m \geq 0$ and $n \geq 1$,

$$\inf_{\|\bar{x}\| \leq C\sqrt{n}} \mathbb{E}^{\bar{x}}(I_{n+[un]})^m \geq \delta b_n^{-dp/2} \mathbb{E}\biggl(\sum_x \prod_{j=1}^p l_j(n, x, \varepsilon)\biggr)^m. \tag{6.4}$$

For each integer $k \geq 0$, let $m \geq 0$ be the integer such that $mp \leq k < p(m+1)$. Applying (6.4) and Hölder's inequality gives

$$\inf_{\|\bar{x}\| \leq C\sqrt{n}} \mathbb{E}^{\bar{x}}(I_{n+[un]})^{(k+p)/p}$$

$$\geq \left[\inf_{\|\bar{x}\|\leq C\sqrt{n}} \mathbb{E}^{\bar{x}}\left(\sum_x \prod_{j=1}^p l_j([un]+n,x)\right)^{m+1}\right]^{(k+p)/((m+1)p)}$$

$$\geq (\delta b_n^{-dp/2})^{(k+p)/((m+1)p)}\left[\mathbb{E}\left(\sum_x \prod_{j=1}^p l_j(n,x,\varepsilon)\right)^{m+1}\right]^{(k+p)/((m+1)p)}$$

$$\geq (\delta b_n^{-dp/2})^{(k+p)/k}\left[\mathbb{E}\left(\sum_x \prod_{j=1}^p l_j(n,x,\varepsilon)\right)^{k/p}\right]^{(k+p)/k}.$$

Therefore,

$$\begin{aligned}
&\left(\frac{b_n}{n}\right)^{((2p-d(p-1))/(2p))(k+p)} \inf_{\|\bar{x}\|\leq C\sqrt{n}} \mathbb{E}^{\bar{x}}(I_{n+[un]})^{(k+p)/p}\\
&\geq \left(\left(\frac{b_n}{n}\right)^{((2p-d(p-1))/(2p))(k+p)} \inf_{\|\bar{x}\|\leq C\sqrt{n}} \mathbb{E}^{\bar{x}}(I_{n+[un]})^{(k+p)/p}\right)^{k/(k+p)} \qquad (6.5)\\
&\geq \delta b_n^{-dp/2}\left(\frac{b_n}{n}\right)^{((2p-d(p-1))/(2p))k}\mathbb{E}\left(\sum_x \prod_{j=1}^p l_j(n,x,\varepsilon)\right)^{k/p},
\end{aligned}$$

where the first step follows from the following rough estimate: As $n$ is sufficiently large,

$$\left(\frac{b_n}{n}\right)^{((2p-d(p-1))/(2p))(k+p)} \inf_{\|\bar{x}\|\leq C\sqrt{n}} \mathbb{E}^{\bar{x}}(I_{n+[un]})^{(k+p)/p} \geq 1, \qquad k=0,1,\dots,$$

which we prove as follows. First, by Jensen's inequality, we need only to show that

$$\left(\frac{b_n}{n}\right)^{(2p-d(p-1))/2} \inf_{\|\bar{x}\|\leq C\sqrt{n}} \mathbb{E}^{\bar{x}}(I_{n+[un]}) \geq 1. \qquad (6.6)$$

Second, similarly to (6.4) (with $m=1$),

$$\begin{aligned}
&\inf_{\|\bar{x}\|\leq C\sqrt{n}} \mathbb{E}^{\bar{x}}(I_{n+[un]})\\
&\quad\geq \delta(\varepsilon)\mathbb{E}\left(\sum_{x\in\mathbb{Z}^d}\prod_{j=1}^p \frac{1}{\#\{B_n'\}}\sum_{k=1}^n \mathbb{1}\left\{\left|\frac{S_j(k)-x}{\sqrt{n}}\right|\leq\varepsilon\right\}\right)\\
&\quad\sim \delta(\varepsilon)n^{(2p-d(p-1))/(2p)}\\
&\qquad\times \mathbb{E}\left(\int_{\mathbb{R}^d}\prod_{j=1}^p \frac{1}{nC_d\varepsilon^d}\sum_{k=1}^n \mathbb{1}\left\{\left|\frac{S_j(k)-[\sqrt{n}x]}{\sqrt{n}}\right|\leq\varepsilon\right\}dx\right) \qquad (n\to\infty),
\end{aligned}$$

where $B_n' = \{x \in \mathbb{Z}^d; |x| \le \varepsilon\sqrt{n}\}$. Therefore, (6.6) follows from the invariance principle which claims

$$\int_{\mathbb{R}^d} \prod_{j=1}^p \frac{1}{nC_d\varepsilon^d} \sum_{k=1}^n \mathbb{1}\left\{\left|\frac{S_j(k) - [\sqrt{n}x]}{\sqrt{n}}\right| \le \varepsilon\right\} dx \xrightarrow{d} \int_{\mathbb{R}^d} \prod_{j=1}^p \widetilde{L}_j(1, x, \varepsilon)\, dx.$$

[Recall that $\widetilde{L}(t, x, \varepsilon)$ is the analogue of $L(t, x, \varepsilon)$ as $W(t)$ is replaced by a Lévy Gaussian process with the same covariance matrix as the random walks.]

Combined with Taylor's expansion, (6.5) implies that for any $0 < \nu < \theta$,

$$\begin{aligned}
&\inf_{\|\bar{x}\| \le C\sqrt{n}} \mathbb{E}^{\bar{x}}\biggl[\left(\frac{b_n}{n}\right)^{(2p-d(p-1))/2} I_{n+[un]} \\
&\qquad\qquad \times \exp\left\{(\theta - \nu)\left(\frac{b_n}{n}\right)^{(2p-d(p-1))/(2p)} (I_{n+[un]})^{1/p}\right\}\biggr] \\
&\quad \ge \delta b_n^{-dp/2} \mathbb{E} \exp\left\{(\theta - \nu)\left(\frac{b_n}{n}\right)^{(2p-d(p-1))/(2p)} \left(\sum_{x \in \mathbb{Z}^d} \prod_{j=1}^p l_j(n, x, \varepsilon)\right)^{1/p}\right\}.
\end{aligned}$$

By Theorem 3.1 and by the fact that $\exp(\theta x^{1/p}) \ge x \exp((\theta - \nu)x^{1/p})$ for large $x$, we have

$$\limsup_{n \to \infty} \frac{1}{b_n} \log \inf_{\|\bar{x}\| \le C\sqrt{n}} \mathbb{E}^{\bar{x}} \exp\left\{\theta\left(\frac{b_n}{n}\right)^{(2p-d(p-1))/(2p)} (I_{n+[un]})^{1/p}\right\} \ge \widetilde{N}_\varepsilon(\theta - \nu).$$

As $\varepsilon \to 0$, the right-hand side approaches

$$\sup_{f \in \mathcal{F}_d} \left\{(\theta - \nu)\left(\int_{\mathbb{R}^d} |f(x)|^{2p}\, dx\right)^{1/p} - \frac{p}{2}\int_{\mathbb{R}^d} \langle \nabla f, \Gamma \nabla f\rangle\, dx\right\}.$$

Replacing $n$ with $k_n \equiv [(1-u)n]$ and noticing that $[uk_n] + k_n \le n$, we have

$$\begin{aligned}
&\limsup_{n \to \infty} \frac{1}{b_n} \log \inf_{\|\bar{x}\| \le C\sqrt{n}} \mathbb{E}^{\bar{x}} \exp\left\{\theta\left(\frac{b_n}{n}\right)^{(2p-d(p-1))/(2p)} I_n^{1/p}\right\} \\
&\quad \ge \sup_{f \in \mathcal{F}_d} \biggl\{(\theta - \nu)(1-u)^{(2p-d(p-1))/(2p)} \left(\int_{\mathbb{R}^d} |f(x)|^{2p}\, dx\right)^{1/p} \\
&\qquad\qquad - \frac{p}{2}\int_{\mathbb{R}^d} \langle \nabla f, \Gamma \nabla f\rangle\, dx\biggr\} \\
&\quad = (1-u)(\theta - \nu)^{2p/(2p-d(p-1))} p^{-(d(p-1))/(2p-d(p-1))} \\
&\qquad \times \det(\Gamma)^{-(p-1)/(2p-d(p-1))} M,
\end{aligned}$$

where the last step follows from the substitution

$$f(x) = \sqrt{|\det(A)|}\, g(Ax)$$

with the $d \times d$ nondegenerate matrix $A$ satisfying

$$A^\tau \Gamma A = (1-u)\left(\frac{\theta-\nu}{p}\right)^{2p/(2p-d(p-1))} \det(\Gamma)^{-(p-1)/(2p-d(p-1))} I_d,$$

where $I_d$ is the $d \times d$ identity matrix. Finally, letting $u \to 0$ and $\nu \to 0$ gives (6.3).

We now prove (6.2) without assuming aperiodicity. Let $0 < \gamma < 1$ be fixed and let $\{\delta_n\}_{n\geq 1}$ be i.i.d. Bernoulli random variables with the common law:

$$\mathbb{P}\{\delta_1 = 0\} = 1 - \mathbb{P}\{\delta_1 = 1\} = \gamma.$$

We assume independence between $\{S(n)\}$ and $\{\delta_n\}$.

Define the renewal sequence $\{\sigma_k\}_{k\geq 0}$ by

$$\sigma_0 = 0 \quad \text{and} \quad \sigma_{k+1} = \inf\{n > \sigma_k; \delta_n = 1\}.$$

Then $\{\sigma_k - \sigma_{k-1}\}_{k\geq 1}$ is an i.i.d. sequence with common distribution

$$\mathbb{P}\{\sigma_1 = n\} = (1-\gamma)\gamma^{n-1}, \qquad n = 1, 2, \dots.$$

Consider the random walk $\{S(\sigma_k)\}_{k\geq 1}$. It is symmetric with covariance

$$\operatorname{Cov}(S(\sigma_1), S(\sigma_1)) = (\mathbb{E}\sigma_1)\Gamma = (1-\gamma)^{-1}\Gamma \tag{6.7}$$

and

$$\mathbb{P}\{S(\sigma_1) = 0\} = (1-\gamma)\sum_{k=1}^{\infty} \gamma^{k-1} \mathbb{P}\{S(k) = 0\} > 0.$$

In particular, $\{S(\sigma_k)\}_{k\geq 1}$ is aperiodic.

Write $\tilde{l}(n,x)$ for the local times of the random walk $\{S(\sigma_k)\}_{k\geq 1}$:

$$\tilde{l}(n,x) = \sum_{k=1}^{n} \mathbb{1}_{\{S(\sigma_k)=x\}}, \qquad x \in \mathbb{Z}^d, n = 1, 2, \dots.$$

Let $\{\delta_n^j\}_{n\geq 1}$ $(1 \leq j \leq p)$ be independent copies of $\{\delta_n\}_{n\geq 1}$ and let the renewal sequences $\{\sigma_k^j\}_{k\geq 0}$ and the local times $\tilde{l}_j(n,x)$ $(1 \leq j \leq p)$ be defined in the obvious way. Write

$$\tilde{I}_n = \sum_{k_1,\dots,k_p=1}^{n} \mathbb{1}_{\{S_1(\sigma_{k_1}^1)=\cdots=S_p(\sigma_{k_p}^p)\}} = \sum_{x\in\mathbb{Z}^d} \prod_{j=1}^{p} \tilde{l}_j(n,x).$$

By what we have proved [i.e., (6.2) under aperiodicity] in the previous step and by (6.7),

$$\begin{aligned} &\liminf_{n\to\infty} \frac{1}{b_n} \log \inf_{\|\bar{x}\|\leq C\sqrt{n}} \mathbb{P}^{\bar{x}}\{\tilde{I}_n \geq \lambda n^{(2p-d(p-1))/2} b_n^{(d(p-1))/2}\} \\ &\qquad \geq -\frac{p}{2}(1-\gamma)^{-1} \det(\Gamma)^{1/d} \kappa(d,p)^{-4p/(d(p-1))} \lambda^{2/(d(p-1))}. \end{aligned} \tag{6.8}$$

On the other hand, notice that

$$\tilde{l}(n,x) = \sum_{k=1}^{\sigma_n} \delta_k \mathbb{1}_{\{S(k)=x\}} \le \sum_{k=1}^{\sigma_n} \mathbb{1}_{\{S(k)=x\}}, \qquad x \in \mathbb{Z}^d, n = 1, 2, \ldots.$$

Consequently, on the events $\{\sigma_n^j < a(1-\gamma)^{-1}n, 1 \le j \le p\}$, where we let $a > 1$, we have

$$\tilde{I}_n \le \sum_{x\in\mathbb{Z}^d} \prod_{j=1}^p l_j(\sigma_n^j, x) \le I_{[a(1-\gamma)^{-1}n]}.$$

So we have, for any $\bar{x} \in (\mathbb{Z}^d)^p$,

$$\begin{aligned}\mathbb{P}^{\bar{x}}&\{\tilde{I}_n \ge \lambda n^{(2p-d(p-1))/2} b_n^{(d(p-1))/2}\}\\ &\le \mathbb{P}^{\bar{x}}\{I_{[a(1-\gamma)^{-1}n]} \ge \lambda n^{(2p-d(p-1))/2} b_n^{(d(p-1))/2}\} + p\mathbb{P}\{\sigma_n \ge \lambda(1-\gamma)^{-1}n\}.\end{aligned}$$

According to Cramér's large deviation [Dembo and Zeitouni (1998), Theorem 2.2.3], there is a $\delta > 0$ such that

$$\mathbb{P}\{\sigma_n \ge a(1-\gamma)^{-1}n\} \le e^{-\delta n}$$

for large $n$. From (6.8), therefore,

$$\begin{aligned}\liminf_{n\to\infty} \frac{1}{b_n} \log \inf_{\|\bar{x}\|\le C\sqrt{n}} \mathbb{P}^{\bar{x}}&\{I_{[a(1-\gamma)^{-1}n]} \ge \lambda n^{(2p-d(p-1))/2} b_n^{(d(p-1))/2}\}\\ &\ge -\frac{p}{2}(1-\gamma)^{-1} \det(\Gamma)^{1/d} \kappa(d,p)^{-4p/(d(p-1))} \lambda^{2/(d(p-1))}.\end{aligned}$$

Replacing $n$ with $[a^{-1}(1-\gamma)n]$ and $\lambda$ with $(a(1-\gamma)^{-1})^{(2p-d(p-1))/2}\lambda$ we obtain

$$\begin{aligned}\liminf_{n\to\infty} \frac{1}{b_n} \log \inf_{\|\bar{x}\|\le C\sqrt{n}} &\mathbb{P}^{\bar{x}}\{I_n \ge \lambda n^{(2p-d(p-1))/2} b_n^{(d(p-1))/2}\}\\ &\ge -\frac{p}{2} a^{(2p-d(p-1))/(d(p-1))} (1-\gamma)^{-2p/(d(p-1))}\\ &\quad\times \det(\Gamma)^{1/d} \kappa(d,p)^{-4p/(d(p-1))} \lambda^{2/(d(p-1))}.\end{aligned}$$

Letting $a \to 1$ and $\gamma \to 0$, we have (6.2). □

REMARK 6.1. Notice that for any $\bar{x} \in (\mathbb{R}^d)^p$ and $t > 0$,

$$\mathcal{L}^{\bar{x}}(\alpha([0,t]^p)) \stackrel{d}{=} \mathcal{L}^{\bar{x}/\sqrt{t}}(t^{(2p-d(p-1))/2}\alpha([0,1]^p)).$$

In particular, (6.1) implies that for any $C>0$ and $\lambda>0$,

$$\liminf_{t\to\infty}\frac{1}{\log\log t} \times \log \inf_{\|\bar{x}\|\le C\sqrt{t}} \mathbb{P}^{\bar{x}}\{\alpha([0,t]^p)\ge \lambda t^{(2p-d(p-1))/2}(\log\log t)^{(d(p-1))/2}\} \ge -\frac{p}{2}\kappa(d,p)^{-4p/(d(p-1))}\lambda^{2/(d(p-1))}. \tag{6.9}$$

This fact is used to prove the lower bound of the law of the iterated logarithm in Theorem 2.3.

**7. The law of the iterated logarithm.** In this section, we prove Theorem 2.3. Due to similarity, we prove only (2.8) in the context of random walks. The proof of the upper bound becomes a standard argument [see, e.g., Chen and Li (2004), Section 6] via the Borel–Cantelli lemma after we take $b_n=\log\log n$ in the moderate deviation given in Theorem 2.2.

To prove the lower bound, let $n_k=k^k$. We first show that for any

$$\lambda<\left(\frac{2}{p}\right)^{(d(p-1))/2}\det(\Gamma)^{-(p-1)/2}\kappa(d,p)^{2p},$$

$$\limsup_{k\to\infty} n_{k+1}^{-(2p-d(p-1))/2}(\log\log n_{k+1})^{-(d(p-1))/2} \times \sum_{x\in\mathbb{Z}^d}\prod_{j=1}^p[l_j(n_{k+1},x)-l_j(n_k,x)]\ge\lambda \qquad \text{a.s.} \tag{7.1}$$

We adapt the notation introduced in Section 6 and consider the $dp$-dimensional random walk $\bar{S}(n)=(S_1(n),\dots,S_p(n))$. By the Markov property and Lévy's Borel–Cantelli lemma [see Breiman (1992), Corollary 5.29], (7.1) holds if we have

$$\sum_k \mathbb{P}^{\bar{S}(n_k)}\{I_{n_{k+1}-n_k}\ge\lambda n_{k+1}^{(2p-d(p-1))/2}(\log\log n_{k+1})^{(d(p-1))/2}\} =\infty \qquad \text{a.s.} \tag{7.2}$$

Indeed, it is easy to see that $\sqrt{n_k\log\log n_k}=o(\sqrt{n_{k+1}-n_k})$ as $k\to\infty$. By the classic Hartman–Wintner law of the iterated logarithm, with probability 1 the events

$$\{\|\bar{S}(n_k)\|\le\sqrt{n_{k+1}-n_k}\},\qquad k=1,2,\dots,$$

eventually hold. Therefore, (7.2) holds if we have

$$\sum_k \inf_{\|\bar{x}\|\le\sqrt{n_{k+1}-n_k}}\mathbb{P}^{\bar{x}}\{I_{n_{k+1}-n_k}\ge\lambda n_{k+1}^{(2p-d(p-1))/2}(\log\log n_{k+1})^{(d(p-1))/2}\}=\infty,$$

which follows from (6.2) in Theorem 6.1 with $b_n = \log\log n$.

Since

$$\sum_{x\in\mathbb{Z}^d}\prod_{j=1}^p l_j(n_{k+1},x) \ge \sum_{x\in\mathbb{Z}^d}\prod_{j=1}^p [l_j(n_{k+1},x) - l_j(n_k,x)],$$

letting

$$\lambda \to \left(\frac{2}{p}\right)^{(d(p-1))/2} \det(\Gamma)^{-(p-1)/2}\kappa(d,p)^{2p}$$

in (7.1) gives the desired lower bound for (2.8).

## APPENDIX

LEMMA A.1. *Under the notation given in Section* 3,

$$\limsup_{m\to\infty}\sup_{g\in\mathcal{F}_d}\left\{\theta\left(\int_{[0,m]^d}\left(\sum_{\mathbf{k}\in\mathbb{Z}^d}(g^2)_\varepsilon(x+m\mathbf{k})\right)^p dx\right)^{1/p} - \frac{1}{2}\int_{\mathbb{R}^d}|\nabla g(x)|^2\,dx\right\} \le M_\varepsilon(\theta). \tag{A.1}$$

PROOF. Let $g \in \mathcal{F}_d$ be fixed and write

$$\bar{g}(x) = \sqrt{\sum_{\mathbf{k}\in\mathbb{Z}^d} g^2(x+m\mathbf{k})}, \qquad x\in\mathbb{R}^d.$$

Then $\bar{g}$ is absolutely continuous and

$$\int_{[0,m]^d}\bar{g}^2(x)\,dx = 1 \quad\text{and}\quad |\nabla\bar{g}(x)|^2 \le \sum_{\mathbf{k}\in\mathbb{Z}^d}|\nabla g(x+m\mathbf{k})|^2. \tag{A.2}$$

Write

$$E = \bigcup_{i=1}^d (\{0\le x_i\le\sqrt{m}+\varepsilon\}\cup\{m-\sqrt{m}-\varepsilon\le x_i\le m\}).$$

By Lemma 3.4 in Donsker and Varadhan (1975), there is an $a\in\mathbb{R}^d$ such that

$$\int_{E^\varepsilon}\bar{g}^2(x+a)\,dx \le \frac{2d(\sqrt{m}+2\varepsilon)}{m},$$

where $E^\varepsilon$ is the $\varepsilon$ neighborhood of $E$. We may assume $a=0$, that is,

$$\int_{E^\varepsilon}\bar{g}^2(x)\,dx \le \frac{2d(\sqrt{m}+2\varepsilon)}{m}, \tag{A.3}$$

for otherwise we may replace $g(\cdot)$ with $g(a+\cdot)$.

Define the function $\phi$ on $\mathbb{R}$ by

$$\phi(\lambda)=\begin{cases}\lambda m^{-1/2}, & 0\le x\le m^{1/2},\\ 1, & m^{-1/2}\le\lambda\le m-m^{1/2},\\ m^{1/2}-\lambda m^{-1/2}, & m-m^{1/2}\le x\le m,\\ 0, & \text{otherwise,}\end{cases}$$

and write

$$\varphi(x)=\phi(x_1)\cdots\phi(x_d),\qquad x=(x_1,\dots,x_d)\in\mathbb{R}^d,$$

$$f(x)=\bar g(x)\varphi(x)\Big/\sqrt{\int_{\mathbb{R}^d}\bar g^2(y)\varphi^2(y)\,dy}=\bar g(x)\varphi(x)/\sqrt{A}\qquad\text{(say).}$$

Then $|\varphi|\le 1$, $|\nabla\varphi|\le\sqrt{d/m}$ and $f\in\mathcal{F}_d$. By (A.2),

$$\begin{aligned}\int_{\mathbb{R}^d}|\nabla f|^2\,dx&=\frac{1}{A}\biggl\{\int_{\mathbb{R}^d}|\nabla\bar g|^2|\varphi|^2\,dx+\int_{\mathbb{R}^d}|\bar g|^2|\nabla\varphi|^2\,dx+2\int_{\mathbb{R}^d}\bar g\varphi\langle\nabla\bar g,\nabla\varphi\rangle\,dx\biggr\}\\ &\le\frac{1}{A}\biggl\{\int_{[0,m]^d}|\nabla\bar g|^2\,dx+\frac{d}{m}\int_{[0,m]^d}|\bar g|^2\,dx\\ &\qquad+2\biggl(\int_{[0,m]^d}|\nabla\bar g|^2|\nabla\varphi|^2\,dx\biggr)^{1/2}\biggr\}\qquad\text{(A.4)}\\ &\le\frac{1}{A}\biggl\{\int_{[0,m]^d}|\nabla\bar g|^2\,dx+\frac{d}{m}+2\sqrt{\frac{d}{m}}\biggl(\int_{[0,m]^d}|\nabla\bar g|^2\,dx\biggr)^{1/2}\biggr\}\\ &\le\frac{1}{A}\biggl\{\biggl(1+\frac{d}{m}\biggr)\int_{[0,m]^d}|\nabla\bar g|^2\,dx+\frac{2d}{m}\biggr\}\\ &\le\frac{1}{A}\biggl\{\biggl(1+\frac{d}{m}\biggr)\int_{\mathbb{R}^d}|\nabla g|^2\,dx+\frac{2d}{m}\biggr\},\end{aligned}$$

where the fourth step follows the inequality $2\theta_1\theta_2\le\theta_1^2+\theta_2^2$.

On the other hand, notice that $(\bar g^2)_\varepsilon(x)=A^{-1/2}(f^2)_\varepsilon(x)$ for all $x\in[0,m]^d\setminus E$:

$$\begin{aligned}&\biggl(\int_{[0,m]^d}\biggl(\sum_{\mathbf{k}\in\mathbb{Z}^d}(g^2)_\varepsilon(x+m\mathbf{k})\biggr)^p dx\biggr)^{1/p}\\ &\qquad=\biggl(\int_{[0,m]^d}|(\bar g^2)_\varepsilon(x)|^p\,dx\biggr)^{1/p}\\ &\qquad\le A\biggl(\int_{[0,m]^d\setminus E}|(f^2)_\varepsilon(x)|^p\,dx\biggr)^{1/p}+\biggl(\int_E|(\bar g^2)_\varepsilon(x)|^p\,dx\biggr)^{1/p}\qquad\text{(A.5)}\end{aligned}$$

$$\le A\biggl(\int_{\mathbb{R}^d}|(f^2)_\varepsilon(x)|^p\,dx\biggr)^{1/p} + \biggl(\frac{2d(\sqrt{m}+2\varepsilon)}{m}\biggr)^{1/p}\sup_x(\bar{g}^2)_\varepsilon(x)^{(p-1)/p},$$

where the last step partially follows from (A.3). By (A.2),

$$(\bar{g}^2)_\varepsilon(x)=\frac{1}{C_d\varepsilon^d}\int_{\{|y-x|\le\varepsilon\}}\bar{g}^2(y)\,dy\le\frac{1}{C_d\varepsilon^d},\qquad x\in\mathbb{R}^d.$$

Combining (A.4) and (A.5) and noticing that $A\le 1$, we obtain

$$\begin{aligned}
&\theta\biggl(\int_{[0,m]^d}\biggl(\sum_{\mathbf{k}\in\mathbb{Z}^d}(g^2)_\varepsilon(x+m\mathbf{k})\biggr)^p dx\biggr)^{1/p}-\frac{1}{2}\int_{\mathbb{R}^d}|\nabla g(x)|^2\,dx\\
&\quad\le A\biggl(1+\frac{d}{m}\biggr)^{-1}\\
&\qquad\times\biggl\{\theta\biggl(1+\frac{d}{m}\biggr)\biggl(\int_{\mathbb{R}^d}|(f^2)_\varepsilon(x)|^p\,dx\biggr)^{1/p}-\frac{1}{2}\biggl(\int_{\mathbb{R}^d}|\nabla f(x)|^2\,dx\biggr)\biggr\}\\
&\qquad+\theta\biggl(\frac{2d(\sqrt{m}+2\varepsilon)}{m}\biggr)^{1/p}\biggl(\frac{1}{C_d\varepsilon^d}\biggr)^{(p-1)/p}+\frac{2d}{m+d}\\
&\quad\le M_\varepsilon\biggl(\theta\biggl(1+\frac{d}{m}\biggr)\biggr)+\theta\biggl(\frac{2d(\sqrt{m}+2\varepsilon)}{m}\biggr)^{1/p}\biggl(\frac{1}{C_d\varepsilon^d}\biggr)^{(p-1)/p}+\frac{2d}{m+d}.
\end{aligned}$$

Taking the supremum on the left-hand side over $g\in\mathcal{F}_d$ and then letting $m\to\infty$ on the both sides, we have (A.1). $\square$

LEMMA A.2. *Under* (1.1),

$$\begin{aligned}
&\sup_{f\in\mathcal{F}_d}\biggl\{\biggl(\int_{\mathbb{R}^d}|f(x)|^{2p}\,dx\biggr)^{1/p}-\frac{1}{2}\int_{\mathbb{R}^d}|\nabla f(x)|^2\,dx\biggr\}\\
&\qquad=\frac{2p-d(p-1)}{2p}\biggl(\frac{d(p-1)}{p}\biggr)^{(d(p-1))/(2p-d(p-1))}\\
&\qquad\quad\times\kappa(d,p)^{4p/(2p-d(p-1))},
\end{aligned}\tag{A.6}$$

*where*

$$\mathcal{F}_d=\biggl\{f\in W^{1,2}(\mathbb{R}^d);\int_{\mathbb{R}^d}|f(x)|^2\,dx=1\biggr\}.$$

PROOF. Write $M$ for the left-hand side of (A.6). For any $f\in\mathcal{F}_d$,

$$\biggl(\int_{\mathbb{R}^d}|f(x)|^{2p}\,dx\biggr)^{1/p}-\frac{1}{2}\int_{\mathbb{R}^d}|\nabla f(x)|^2\,dx$$

$$\leq \kappa(d,p)^2 |\nabla f|_2^{(d(p-1))/p} - \frac{1}{2}|\nabla f|_2^2 \leq \sup_{\theta>0}\left\{\kappa(d,p)^2\theta^{(d(p-1))/p} - \frac{1}{2}\theta^2\right\}$$
$$= \frac{2p-d(p-1)}{2p}\left(\frac{d(p-1)}{p}\right)^{(d(p-1))/(2p-d(p-1))}\kappa(d,p)^{4p/(2p-d(p-1))}.$$

Hence

$$M \leq \frac{2p-d(p-1)}{2p}\left(\frac{d(p-1)}{p}\right)^{(d(p-1))/(2p-d(p-1))}\kappa(d,p)^{4p/(2p-d(p-1))}.$$

On the other hand, for any $C < \kappa(d,p)$, there is a $g$ such that

$$\|g\|_{2p} > C\|\nabla g\|_2^{(d(p-1))/(2p)} \cdot \|g\|_2^{1-(d(p-1))/(2p)}.$$

By homogeneity, we may assume $\|g\|_2 = 1$. Given $\lambda > 0$, let $f(x) = \lambda^{d/2}g(\lambda x)$. Then $\|f\|_2 = 1$, $\|\nabla f\|_2 = \lambda\|\nabla g\|_2$ and

$$\|f\|_{2p} = \lambda^{(d(p-1))/(2p)}\|g\|_{2p} > C(\lambda\|\nabla g\|_2)^{(d(p-1))/(2p)}.$$

Hence,

$$M \geq \|f\|_{2p}^2 - \tfrac{1}{2}\|\nabla f\|_2^2 > C^2(\lambda\|\nabla g\|_2)^{(d(p-1))/p} - \tfrac{1}{2}(\lambda\|\nabla g\|_2)^2.$$

Since $\lambda > 0$ can be arbitrary,

$$M \geq \sup_{\theta>0}\left\{C^2\theta^{(d(p-1))/p} - \frac{1}{2}\theta^2\right\}$$
$$= \frac{2p-d(p-1)}{2p}\left(\frac{d(p-1)}{p}\right)^{(d(p-1))/(2p-d(p-1))}C^{4p/(2p-d(p-1))}.$$

Letting $C \to \kappa(d,p)$ gives the desired lower bound. □

**Acknowledgments.** I thank R. F. Bass, A. Dorogovtsev, D. Khoshnevisan, S. Kwapien, W. Li and J. Rosen for their interest in this work and for their helpful comments. Also, I am grateful for many good suggestions from the referee.

## REFERENCES

Bass, R. F. and Chen, X. (2004). Self-intersection local time: Critical exponent and laws of the iterated logarithm. *Ann. Probab.* **32** 3221–3247. MR2094444

Bass, R. F. and Khoshnevisan, D. (1992). Local times on curves and uniform invariance principles. *Probab. Theory Related Fields* **92** 465–492. MR1169015

Bass, R. F. and Khoshnevisan, D. (1993). Intersection local times and Tanaka formulas. *Ann. Inst. H. Poincaré Probab. Statist.* **29** 419–451. MR1246641

Breiman, L. (1992). *Probability*. SIAM, Philadelphia. MR1163370

Carlen, E. A. and Loss, M. (1993). Sharp constant in Nash's inequality. *Internat. Math. Res. Notices* **1993** 213–215. MR1230297

Chen, X. and Li, W. (2004). Large and moderate deviations for intersection local times. *Probab. Theory Related Fields* **128** 213–254. MR2031226

Cordero-Erausquin, D., Nazaret, B. and Villani, C. (2004). A mass-transportation approach to sharp Sobolev and Gagliardo–Nirenberg inequalities. *Adv. Math.* **182** 307–332. MR2032031

Darling, D. A. and Kac, M. (1957). On occupation times for Markoff processes. *Trans. Amer. Math. Soc.* **84** 444–458. MR84222

de Acosta, A. (1985). Upper bounds for large deviations of dependent random vectors. *Z. Wahrsch. Verw. Gebiete* **69** 551–565. MR791911

Del Pino, M. and Dolbeault, J. (2002). Best constants for Gagliardo–Nirenberg inequalities and applications to nonlinear diffusions. *J. Math. Pures Appl.* **81** 847–875. MR1940370

Del Pino, M. and Dolbeault, J. (2003). The optimal Euclidean $L^p$-Sobolev logarithmic inequality. *J. Funct. Anal.* **197** 151–161. MR1957678

Dembo, A. and Zeitouni, O. (1998). *Large Deviations Techniques and Applications*, 2nd ed. Springer, New York. MR1619036

Donsker, M. D. and Varadhan, S. R. S. (1974). Asymptotic evaluation of certain Wiener integrals for large time. In *Functional Integration and Its Applications* (A. M. Arthurs, ed.) 15–33. Clarendon Press, Oxford. MR486395

Donsker, M. D. and Varadhan, S. R. S. (1975). Asymptotics for the Wiener sausage. *Comm. Pure Appl. Math.* **28** 525–565. MR397901

Dunford, N. and Schwartz, J. T. (1988). *Linear Operators. Part I. General Theory.* Wiley, New York. MR1009162

Dvoretzky, A., Erdös, P. and Kakutani, S. (1950). Double points of paths of Brownian motions in $n$-space. *Acta Sci. Math. (Szeged)* **12** 75–81. MR34972

Dvoretzky, A., Erdös, P. and Kakutani, S. (1954). Multiple points of paths of Brownian motion in the plane. *Bulletin of the Research Council of Israel* **3** 364–371. MR67402

Fitzsimmons, P. J. and Pitman, J. (1999). Kac's moment formula and the Feymann–Kac formula for additive functionals of a Markov process. *Stochastic Process. Appl.* **79** 117–134. MR1670526

Geman, D., Horowitz, J. and Rosen, J. (1984). A local time analysis of intersections of Brownian paths in the plane. *Ann. Probab.* **12** 86–107. MR723731

Khoshnevisan, D. (2003). Intersections of Brownian motions. *Exposition. Math.* **21** 97–114. MR1978059

Khoshnevisan, D., Xiao, Y. and Zhong, Y. (2003a). Local times of additive Lévy processes. *Stochastic. Process. Appl.* **104** 193–216. MR1961619

Khoshnevisan, D., Xiao, Y. and Zhong, Y. (2003b). Measuring the range of an additive Lévy process. *Ann. Probab.* **31** 1097–1141. MR1964960

König, W. and Mörters, P. (2002). Brownian intersection local times: Upper tail asymptotics and thick points. *Ann. Probab.* **30** 1605–1656. MR1944002

Le Gall, J.-F. (1986a). Propriétés d'intersection des marches aléatoires. I. Convergence vers le temps local d'intersection. *Comm. Math. Phys.* **104** 471–507. MR840748

Le Gall, J.-F. (1986b). Propriétés d'intersection des marches aléatoires. II. Étude des cas critiques. *Comm. Math. Phys.* **104** 509–528. MR840749

Le Gall, J.-F. (1990). Some properties of planar Brownian motion. *Lecture Notes in Math.* **1527** 111–235. Springer, Berlin. MR1229519

Le Gall, J.-F. (1994). Exponential moments for the renormalized self-intersection local time of planar Brownian motion. *Séminaire de Probabilités XXVIII. Lecture Notes in Math.* **1583** 172–180. Springer, Berlin. MR1329112

LE GALL, J.-F. and ROSEN, J. (1991). The range of stable random walks. *Ann. Probab.* **19** 650–705. MR1106281

LEVINE, H. A. (1980). An estimate for the best constant in a Sobolev inequality involving three integral norms. *Ann. Mat. Pura Appl.* **124** 181–197. MR591555

MANSMANN, U. (1991). The free energy of the Dirac polaron, an explicit solution. *Stochastics Stochastics Rep.* **34** 93–125. MR1104424

MARCUS, M. B. and ROSEN, J. (1997). Laws of the iterated logarithm for intersections of random walks on $Z^4$. *Ann. Inst. H. Poincaré Probab. Statist.* **33** 37–63. MR1440255

REMILLARD, B. (2000). Large deviations estimates for occupation time integrals of Brownian motion. In *Stochastic Models* 375–398. Amer. Math. Soc., Providence, RI. MR1765021

ROSEN, J. (1990). Random walks and intersection local time. *Ann. Probab.* **18** 959–977. MR1062054

ROSEN, J. (1997). Laws of the iterated logarithm for triple intersections of three-dimensional random walks. *Electron. J. Probab.* **2** 1–32. MR1444245

WEINSTEIN, M. I. (1983). Nonlinear Schrödinger equations and sharp interpolation estimates. *Comm. Math. Phys.* **87** 567–576. MR691044

ZIEMER, W. P. (1989). *Weakly Differentiable Functions*. Springer, New York. MR1014685

DEPARTMENT OF MATHEMATICS
UNIVERSITY OF TENNESSEE
KNOXVILLE, TENNESSEE 37996-1300
USA
E-MAIL: xchen@math.utk.edu
URL: www.math.utk.edu/~xchen